\documentclass[journal]{IEEEtran}
\usepackage{color}
\usepackage{amssymb}
\usepackage{amsmath}
\usepackage{yhmath}
\usepackage{ulem} 
\usepackage{graphicx}
\usepackage{multirow}
\usepackage{multicol}
\usepackage{mdwlist}
\usepackage[colorlinks,linkcolor=blue]{hyperref}
\usepackage{booktabs}
 \usepackage{ulem}

\usepackage{environ}
\NewEnviron{myequation}{%
\begin{equation}
\scalebox{.8}{$\BODY$}
\end{equation}
}

\begin{document}

\title{Towards Fairness-Aware Multi-Objective Optimization}

\author{Guo~Yu,~\IEEEmembership{Member,~IEEE,}         
		Lianbo~Ma,~\IEEEmembership{Member,~IEEE,}
		Wei Du,       	\\              
      Wenli~Du, 
      Yaochu~Jin,~\IEEEmembership{Fellow,~IEEE}

\thanks{Guo Yu, Wei Du, and  Wenli Du are with the Key Laboratory of Smart Manufacturing in Energy Chemical Process, Ministry of Education, East China University of Science and Technology, Shanghai 200237, China.  E-mail: (guoyu;~duwei0203;~wldu)@ecust.edu.cn

Lianbo Ma is with Software College, Northeastern University, Shenyang, 110819, China. E-mail:malb@swc.neu.edu.cn

Yaochu Jin is with the Faculty of Technology, Bielefeld University, 33619 Bielefeld, Germany. He is also with the Department of Computer Science, University of Surrey, Guildford, Surrey GU2 7XH, UK. E-mail: yaochu.jin@uni-bielefeld.de. (\textit{Corresponding authors: Yaochu Jin, Wenli Du})
}
}

\maketitle
\begin{abstract}
Recent years have seen the rapid development of fairness-aware machine learning in mitigating unfairness or discrimination in decision-making in a wide range of applications. However, much less attention has been paid to the fairness-aware multi-objective optimization, which is indeed commonly seen in real life, such as fair resource allocation problems and data-driven multi-objective optimization problems. This paper aims to illuminate and broaden our understanding of multi-objective optimization from the perspective of fairness. To this end, we start with a discussion of user preferences in multi-objective optimization and then explore its relationship to fairness in machine learning and multi-objective optimization. Following the above discussions, representative cases of fairness-aware multi-objective optimization are presented, further elaborating the importance of fairness in traditional multi-objective optimization, data-driven optimization and federated optimization. Finally, challenges and opportunities in fairness-aware multi-objective optimization are addressed. We hope that this article makes a small step forward towards understanding fairness in the context of optimization and promote research interest in fairness-aware multi-objective optimization.

\end{abstract}


\IEEEpeerreviewmaketitle

%
\IEEEpeerreviewmaketitle

\section{Introduction}\label{Sectionofintro}
\IEEEPARstart{R}{cent} decades have seen the rapid development of fairness-aware machine learning (FML) across various areas \cite{pessach2022review}, such as community detection \cite{9073602}, graph or word embedding \cite{WEnips2016,WD2019}, classification \cite{zafar2017fairness}, regression \cite{berk2017convex}, and recommendation systems \cite{yang2017measuring,biega2018equity}. The reason lies in the fact that loads of evidence have indicated machine learning (ML) impacts the interests of unprivileged people when only accuracy of ML models is required in decision making \cite{shi2021survey,angwin2016machine,jabbari2017fairness}. On the other hand, an outpouring of research has also shown that biased search or models can lead to the manipulation of decision-making in election \cite{epstein2015search} and users' perception of things \cite{braaten2011trust,babaei2021analyzing}, an unfair distribution of opportunities (or resources) on healthcare \cite{ledford2019millions}, reinforcement of social stereotypes \cite{kay2015unequal}, discrimination or misclassification based on sensitive attributes in facial recognition \cite{hao2019making}, and criminal justice \cite{berk2021fairness}. We can conclude that unfairness from ML models may produce negative effect on every aspect of our lives. As a consequence, ML is required to not only excel at aiding users in decision making but also bear the social responsibility of being fair \cite{GaoSh19,baeza2018bias}. 

When fairness is taken into account, an alarming issue of ML is the dilemma between model accuracy and fairness, since the fairness may impact the accuracy \cite{speicher2018unified,fairevo2021}. In addition, some fairness measures are conflicting with each other. For example, the increase of the individual fairness may lead to the decrease of the group fairness \cite{chouldechova2017fair}. 

\begin{figure}[t]
\centering
  \includegraphics[width=.75\linewidth]{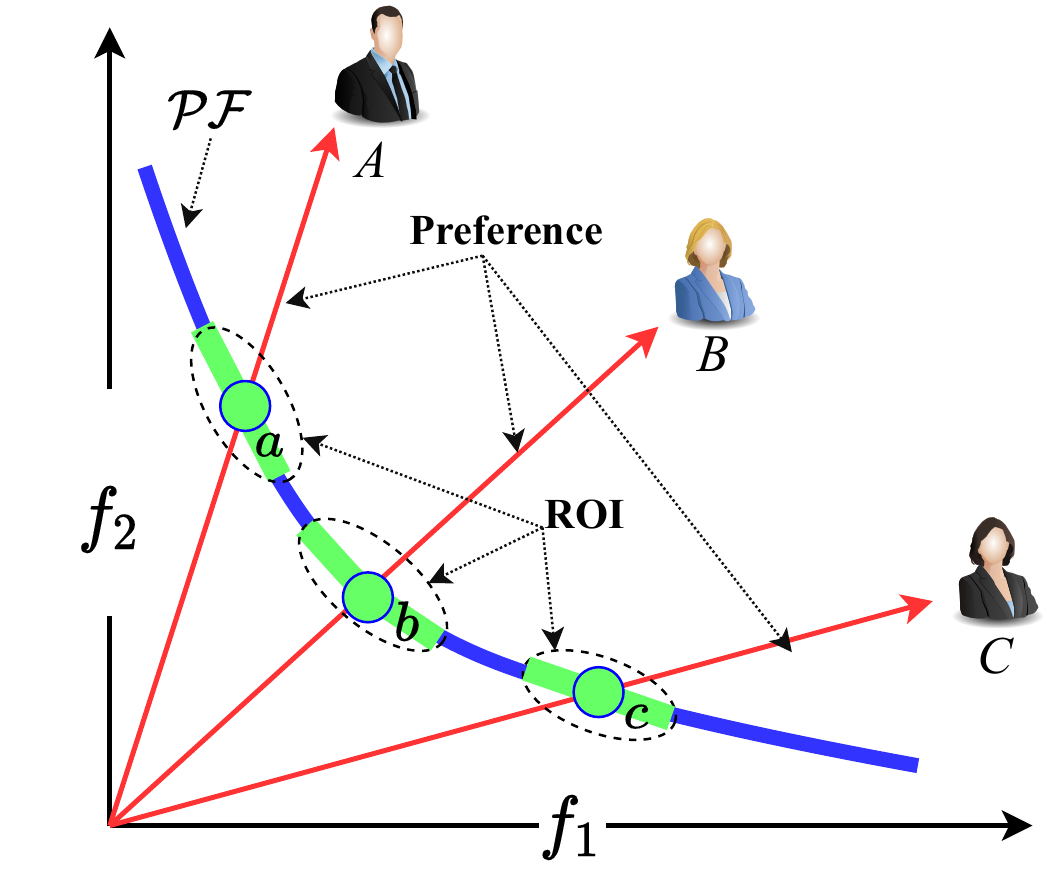}\\
  \caption{An illustration of PMO where preferences based on weight vectors are specified by the DM ($A$, $B$, $C$). ROIs are the regions of interest, and solutions $a$, $b$ and $c$ are the corresponding best desired solution in each of these ROIs.}\label{plot-prefer}
\end{figure}

Interestingly, we have found that FML and preference-driven multi-objective optimization (PMO) \cite{miettinen2012nonlinear,minireview2017,GuoSurvey2021} share many similarities. They both involve human in the loop, aiming to satisfy the users or decision-makers (DMs) with fair and accurate model or solutions satisfying the DM's preferences. In PMO, the preferences from DMs are usually applied in the multi-objective optimization  \cite{minireview2017} to guide the search towards the regions of interest (ROIs) if DMs are able to articulate their preferences \cite{8790106}.  From the perspective of PMO, the best trade-off solution between the utilities (e.g., accuracy and model complexity) and fairness metrics in ML can be taken as the best desired solution in the ROI \cite{adra2007comparative} of the Pareto front. In Fig. \ref{plot-prefer}, an example of the PMO is presented, where solutions ($a$, $b$ and $c$) are the best desired solutions in the corresponding ROIs under different reference vectors. Therefore, when the objectives ($f_1$ and $f_2$) are set to the utility and fairness of ML, then the weight vectors can be extracted according to the relative importance between the objectives, and the solutions are the best trade-off under different weights \cite{fairevo2021,zhang2020joint,valdivia2021fair}. In a joint optimization between fairness and utilities of ML, Zhang et al. \cite{zhang2020joint} argue that ensuring fairness in decision making is not just a matter of creating new algorithms that are ``fairer'' but it is an iterative and interactive process that involves multiple stakeholders. Therefore, integrating preferences from the users into the fairness-aware optimization may alleviate the issue of unfairness in ML \cite{ruchte2021scalable}. 

The notions of fairness in FML are considered in various domains of applications such as classification, clustering, regression (or prediction), and ranking  \cite{chouldechova2018frontiers,mehrabi2021survey,barocas-hardt-narayanan, zehlike2021fairness}. However, the notions of fairness may also appear in different optimization scenarios, on the condition that the biased optimization disproportionately violates the expectations of  certain groups of people such as their interests or preferences. Hence, fairness-aware multi-objective optimization (FMO) is of paramount significance and therefore we believe FMO will be an increasingly important research topic.

For example, in the optimization of dynamic prices for electric vehicle charging, Limmer and Dietrich \cite{limmer2018optimization} designed a FMO framework by considering the objectives of minimization of unfairness of the dynamic prices, maximization of the expected profit, and the maximization of price. The unfairness comes from that it is perceived as unfair to the people who experience different prices for similar products \cite{2004TheUnfairPrice}. Their experimental results show that the increase of the fairness does not impact the expected profit to some extent compared with single optimization of the profit. That is, this example has verified that it is possible to maximize the profit of the charging station operator while mitigating the unfairness issue in dynamic prices for the customers. FMO is an emerging topic in the optimization of dynamic prices for electric vehicle, and more studies are deserved to the problem by considering extra objectives such as the minimization of the number of declining customers and rejection rate \cite{8285280}. 

In addition, with the development of Cloud computing and Internet of Things, the resource allocation with respect to the fairness has attracted many attentions \cite{FairCAS2007}. In the resource allocation problem, the resource providers bid to offer their services in a reverse auction-based cloud market, from which cloud users purchase the resources. Although the users are able to benefit from this type of auction to get their resources at a minimum price, it may result in only several cloud providers dominating the whole cloud market. The reason lies in the fact that most customers show more trust in existing providers. In this situation, the fairness embodies in giving more chances to the ``always losing'' providers \cite{fairRA2021}. Another fairness defined in \cite{ghodsi2011dominant} is to make sure that the users have the equalized maximum ratio of any resource they have in the resource pool, which is also extended to solve the multi-resource allocation in heterogeneous cloud computing systems \cite{6919321}. Although some studies \cite{6211077,JIANG2016239} have taken the cloud resource allocation problem as a multi-objective optimization problem (MOP), the notions of fairness are not taken into account.

In recent years, fairness-aware federated learning has become a hot topic in ML \cite{shi2021survey} for the following reasons. Firstly, the participants or clients are self-interested, and they are not willing to share their data without incentive rewards. Secondly, fairness is required for a sustainable development of the federated learning ecosystem where there is no discrimination against any individuals or groups with specific attributes. Different types of definitions for fairness have been proposed. For example, performance fairness is defined so that more emphasis is put on the uniform accuracy distribution across the participants, instead of only considering the interest of the central controller. Collaboration fairness aims to ensure that participants with  larger contributions will receive higher rewards. Model fairness shows that there is no discrimination of the model against specific groups of people. Meanwhile, the above notions of fairness also exist in federated optimization \cite{xu2021federatedEA,xu2021federatedMOEA}, where data subject to privacy restrictions is collected in a distributed way for data-driven optimization. So far, few studies have ever investigated the fairness-aware federated optimization. 

To the best of our knowledge, no comprehensive discussions about FMO have been reported in the literature. To fill the gap, this survey paper discusses fairness in multi-objective optimization together with FML and PMO to reveal the connections and differences between these three research topics. The main contributions of this paper are as summarized as follows.
\begin{enumerate}
  \item For the first time, the importance of fairness in multi-objective optimization is pointed out and explicitly emphasized, and the relationship between preferences and fairness in multi-objective optimization is explored.
  \item A number of representative case studies of FMO are highlighted to gain deeper insights into fairness in multi-objective optimization, data-driven optimization, and federated optimization.
  \item Challenges and opportunities in FMO and fairness-aware multi-objective ML are identified and promising future topics are suggested.
\end{enumerate}

The rest of the paper is organized as follows. Section \ref{Sectionofmoti} briefly introduces the concept of fairness in multi-objective optimization. Then, Section \ref{SectionofPMO} presents the preliminaries preference-driven multi-objective optimization, followed by a discussion of the relationship between preference and fairness in Section \ref{SectionofP2F}. Section \ref{SectionofFMO} elaborates the notions of fairness in various multi-objective optimization frameworks. Then a number of challenges and promising research opportunities are presented in Section \ref{SectionofCO}. Finally, Section \ref{SectionofConclusion} concludes the paper.

\begin{figure*}[t]
		\centering
		\begin{tabular}{@{}c@{}c@{}}
			\includegraphics[width=.45\linewidth]{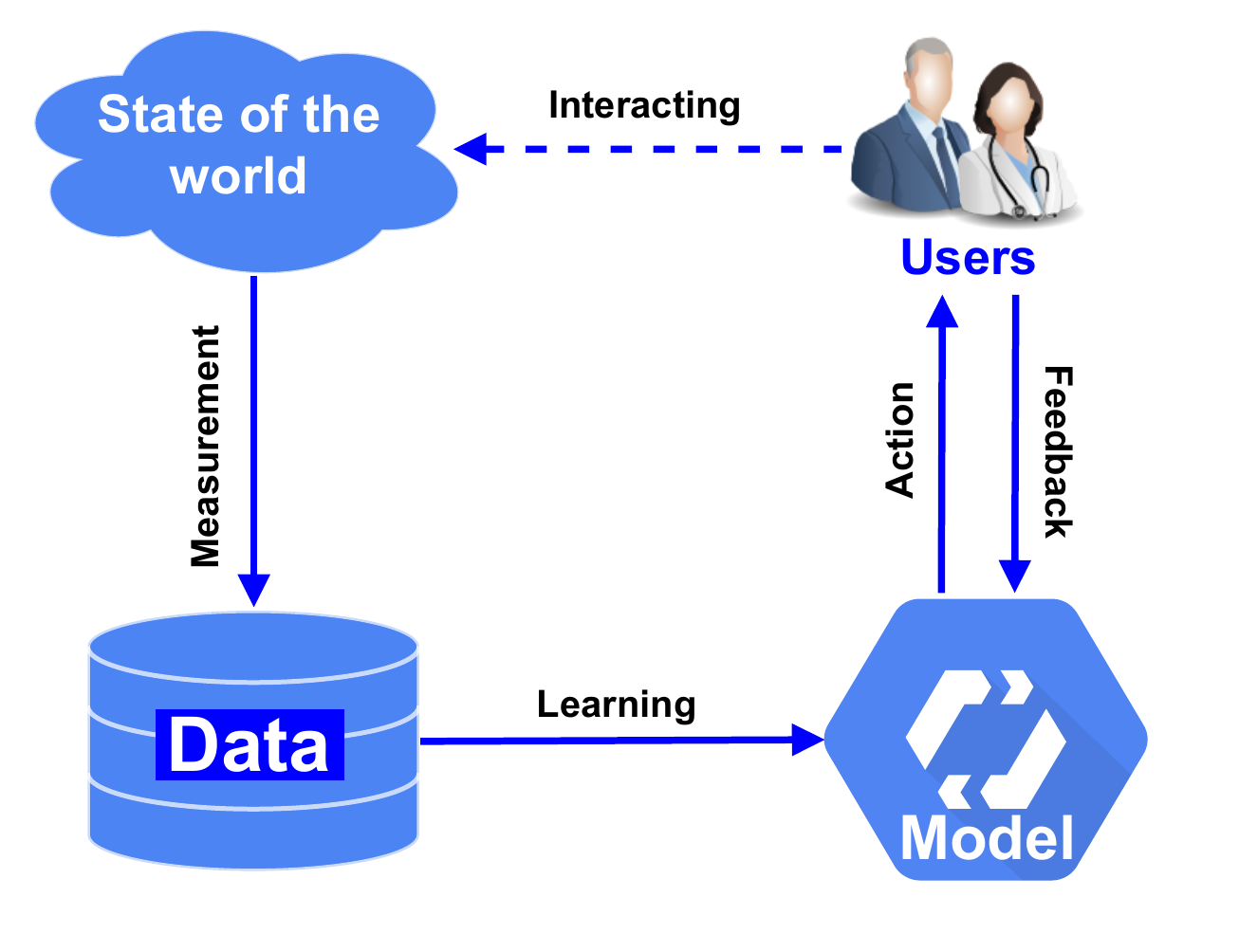}&\includegraphics[width=.42\linewidth]{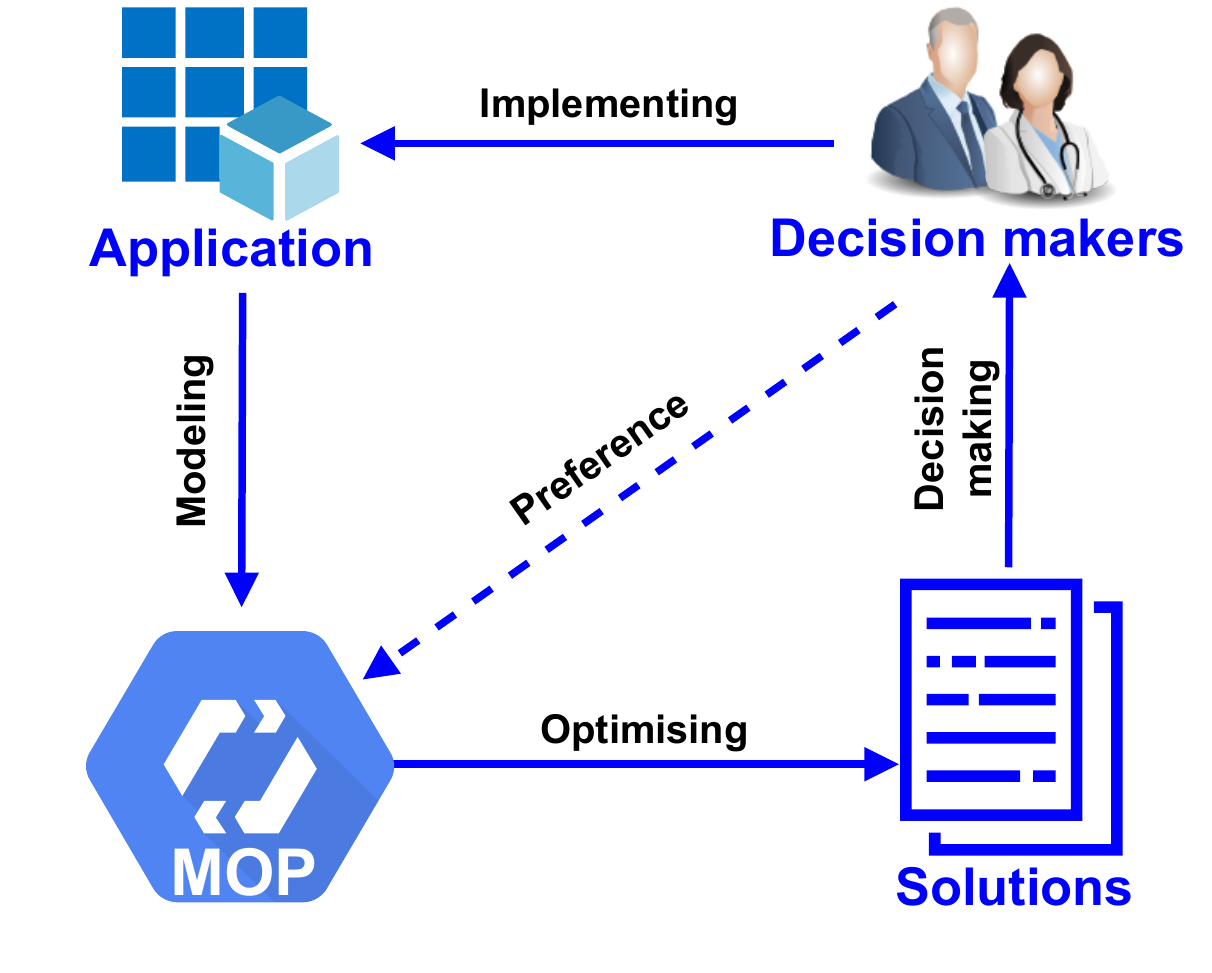}\\
			(a) The loop of ML        & 			(b) The loop of PMO
		\end{tabular}
		\caption{An illustration of the pipelines of ML and PMO. }\label{plot-MLMOP}
\end{figure*}

\section{Motivation}\label{Sectionofmoti}
To discuss fairness in PMO in relation to FML, we plot the pipelines of ML and PMO in Fig. \ref{plot-MLMOP}. From Fig. \ref{plot-MLMOP} (a), we can see that the data used for building ML model is abstracted from the state-of-the-world. After that, the ML model is trained on the data to learn the knowledge in the data. For example, a convolutional neural network can be trained for face recognition \cite{hu2015face} on a set of face images with labels. Once the model is trained, it can then be used to make inferences (predictions), and further refined via the feedback from the users. In some cases, the users may interact with the world through the ML model such as a recommendation system \cite{davidson2010youtube}. In PMO, as shown in Fig. \ref{plot-MLMOP} (b), an MOP is extracted from a real-world application either by mathematical modelling or data-driven surrogate modelling \cite{8456559,jin2021data}. After that, the MOP is optimized by a multi-objective optimization method and consequently a set of representative solutions will be obtained \cite{miettinen2012nonlinear}. Next, the DM will select a small subset of preferred solutions from the obtained solutions \cite{minireview2017,GuoSurvey2021}. The DMs may also interact with the optimization process by incorporating their feedback (e.g., preferences) into the search process to tune the optimization model \cite{zheng2017decomposition,8412189}. Finally, the most preferred solutions will be implemented on the real application. From the above introduction, we can see that the two pipelines share many similarities. For example, both pipelines are interacting with human users to satisfy the users or DMs with an accurate model or preferred solutions. In addition, the model or solutions are obtained through training the model on the data or optimizing an MOP. 

Although there are many similarities between both pipelines, much more attention has been paid to FML \cite{shi2021survey,angwin2016machine,jabbari2017fairness} than FMO.  On the other hand, FMO has become increasingly important in many optimization problems such as the optimization of dynamic prices for electric vehicle charging \cite{limmer2018optimization} and resources allocation optimization \cite{FairCAS2007}. Many concerns about unfairness in multi-objective optimization have arisen, however, little has been done to address these concerns. For example, different distributions of training data may result in surrogates of an MOP with different levels of fidelity \cite{jin2021data}, which may produce biases in the performance and result in biased search towards uninterested regions. In federated optimization \cite{xu2021federatedMOEA} where the data subject to privacy can be collected only in a distributed way, performance fairness \cite{shi2021survey}, collaboration fairness \cite{shi2021survey}, and model fairness \cite{shi2021survey} will also become an important issue. In addition, optimization processes with/out considering fairness will result in different solutions \cite{2004TheUnfairPrice}. A decision making stage without considering fairness may also lead to discrimination issues \cite{8955593}.

To conclude, fairness is of great importance in both ML and multi-objective optimization, and biased ML models or solutions will violate the expectations and interest of certain groups of users or DMs. Thus, fairness in multi-objective optimization will be increasingly important, just as in FML. 

\section{Preliminaries}\label{SectionofPMO}
In this section, some background of the traditional MO and PMO \cite{minireview2017,GuoSurvey2021,8710313} is presented before we introduce fairness-aware optimization.

\subsection{Multi-objective optimization}
MOPs involving multiple conflicting objectives are ubiquitous in the real world, to which a set of trade-off solutions will be found \cite{miettinen2012nonlinear}. Without loss of generality, an MOP can be formulated as follows:
\begin{equation}\label{MOP}
\left\{ \begin{array}{ll}
   \textrm{\it{minimize} }    & \mathbf{F}({\bf{x}})=(f_1({\bf{x}}),\cdots,f_m({\bf{x}}))^{T},\\
   \textrm{\it{subject to} }  & g_j({\bf{x}})\leq 0, j=1,\ldots,J,\\
								& h_k({\bf{x}})= 0, k=1,\ldots,K,\\		  
								&\quad {\bf{x}}\in \Omega,
\end{array} \right.
\end{equation}
where ${\bf{x}}=(x_1,\cdots,x_n)$ is the decision vector. $\Omega \subseteq \mathbb{R}^n$ is the decision space, and $n$ is the number of decision variables. $\mathbf{F}:\Omega \rightarrow \mathbb{Y}$ consists of $m$ objectives. $\mathbb{Y} \subseteq \mathbb{R}^m$ is the objective space. $g$ and $h$ are the inequality constraints and equality constraints, respectively. $J$ and $K$ are the corresponding number of constraints. In multi-objective optimization, a solution $\bf{x_1}$ is said to Pareto dominate another solution $\bf{x_2}$ (denoted by $\bf{x_1}\prec \bf{x_2}$), if and only if $\forall i \in \{1,\ldots,m\}$, $f_i({\bf{x_1}})\leq f_i(\bf{x_2})$ and $\exists j\in \{1,\ldots,m\}$, $f_j({\bf{x_1}}) < f_j(\bf{x_2})$. The Pareto optimal set ($\mathcal{PS}$) of an MOP is formed by all Pareto optimal solutions, i.e., $\mathcal{PS} = \{\bf{x}\in \Omega|\nexists \bf{y}\in \Omega, y\prec \bf{x}\}$, where a solution $\bf{x}^*$ is Pareto optimal if there is no solution $\bf{x}\in \Omega$ that dominates $\bf{x}^*$. The mapping of the PS in the objective space is called the Pareto optimal front ($\mathcal{PF}$), i.e., $\mathcal{PF} = \{\mathbf{F}({\bf{x}})|{\bf{x}}\in \mathcal{PS}\}$.  

In traditional multi-objective optimization, a set of trade-off solutions of the $\mathcal{PF}$ is commonly required to be obtained if no specific preference is given \cite{miettinen2012nonlinear}. If the DM has \textit{a priori} knowledge of the optimization problem, e.g., the relative importance about the different objectives, such knowledge can be articulated as preference information to be incorporated into the optimization process \cite{2011Review,BECHIKH2015141}, leading to the identification of the solutions in the ROIs \cite{minireview2017}. Therefore, the preference-driven methods in evolutionary multi-objective optimization aim to get a representative subset of the Pareto optimal solutions in the ROI \cite{J2008Multiobjective}. By contrast, the methods in multi-criteria decision making community commonly use mathematical programming methodologies to find one most preferred Pareto optimal solution \cite{miettinen2012nonlinear}. 

\subsection{Preference-driven multi-objective optimization}
\subsubsection{\textit{A priori}, interactive and \it{a posteriori} multi-objective optimization}
According to the time when the DM interact with the optimization, i.e., before, during, or after the optimization, the preference-driven methods can be classified into {\it{a priori}}, {\it{interactive}}, {\it{a posteriori}} approaches \cite{miettinen2012nonlinear}. 
\begin{itemize}
\item In the {\it{a priori}} optimization, the DM interact with the optimization by using his/her preferences to guide the search before the optimization starts. An underlying hypothesis is that the DM has sufficient \textit{a priori} knowledge, and is able to accurately articulate the preferences \cite{minireview2017}.

\item The {\it{interactive}} optimization enables the DM to interactively guide the optimization by tuning the DM's preferences on the basis of the already obtained solutions \cite{8710313,8412189}. With the guidance of the progressively provided preferences, the interactive multi-objective optimization will gradually find the preferred solutions. The main difficulty of the interactive methods lies in the fact that the tuning process is arduous and sometimes intractable \cite{8790106,8412189}.  

\item The {\it{a posteriori}} optimization is performed after a representative solution set of the Pareto front ($\mathcal{PF}$) is obtained by any optimization methods \cite{YU2017689,9332241,9484680}. The main idea of the \textit{a posteriori} optimization is to choose a small number of the solutions from the provided solution set according to the DM's preference \cite{miettinen2012nonlinear}. The \textit{a posteriori} optimization will become increasingly difficult when the number of objectives increases, since it is resource-intensive and time-consuming to get a good representative solution set of the $\mathcal{PF}$ in a high-dimensional objective space \cite{8790106,2003Performance,ishibuchi2008evolutionary}. 
\end{itemize}

\subsubsection{User preference modelling}
Preference models are the links or interfaces between the user preferences and the optimization methods. Following the study \cite{8412189}, we categorize the commonly seen preference models into value function (or utility function), dominance relation, and decision rules.

\begin{itemize}
\item {\it{Value function:}} Value functions are commonly seen as scalar (or aggregation) functions of all objectives, which transform the comparison of all objectives into the ranking of the fitness values during the environmental selection \cite{2002Axiomatization}. Widely used scalarizing functions include the weighted sum approach \cite{miettinen2012nonlinear}, Tchebycheff approach \cite{miettinen2012nonlinear}, penalty-based intersection approach \cite{4358754,zheng2017decomposition}, angle penalized distance \cite{7386636}, weighted Euclidean distance approach \cite{deb2006reference}, achievement scalarizing functions \cite{knee9209056,ruiz2008additive}, expected marginal distance functions \cite{Bhattacharjee2017Bridging}, and the minimum Manhattan distance approach \cite{7465803,8933061}. The parameters embedded in these functions are given by the DM to indicate their preferences in the models. 
\item {\it{Dominance relation:}} Dominance relations, many of which are variants of the Pareto dominance, can be used to describe the relationship between the solutions to reflect the DM's preferences, such as the $g$-dominance \cite{2009gdominance},  polyhedral cones based dominance \cite{cone2010}, $r$-dominance \cite{2010rdominance}, $ar$-dominance \cite{8552449}, and localized dominance \cite{G3321930,9139367}. Similarly, the parameters in the dominance relations represent the DM's preferences of solutions or ROIs. 
\item {\it{Decision rules:}} Decision rules aim to model the DM's pairwise solution comparisons so that the solutions can be differentiated according to the defined preference rules or relations \cite{lopez2009study}, such as the rough set based approach \cite{2001Rough,5585982,Greco2008}, fuzzy rules \cite{yaochu2002fuzzy,2012Expressing}, and preference rules \cite{2011Relational,2010Brain}. 
\end{itemize}

In summary, preference models bridge the gap between the DM and the optimization algorithm in various ways, in order to meet different situations and preference articulations. 

\begin{figure}[t]
\centering
  \includegraphics[width=.65\linewidth]{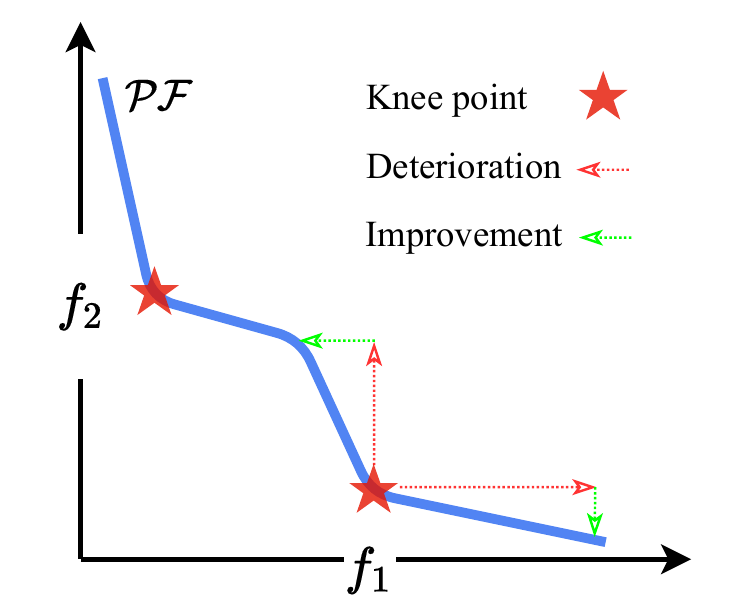}\\
  \caption{An example of the knee points of $\mathcal{PF}$.}\label{plot-knees}
\end{figure}

\subsubsection{User preference articulations}
According to whether a DM can properly articulate his or her preference, the preference information can be classified into explicit (or informed) and implicit (or uninformed) preference \cite{minireview2017,GuoSurvey2021}. 
  
{\it{{Explicit preference:}}} The explicit preference information can be divided into three categories according to whether the comparisons are performed on objectives or solutions \cite{8412189}.

\begin{itemize}
\item {\it{Expectation}}: Expectation is commonly embodied by a goal or a reference point that reflects the aspiration level on each objective the DM intends to achieve \cite{deb2006reference}. The DM can freely set the reference point which may be feasible or infeasible, based on the \textit{a priori} knowledge of the ranges of objectives \cite{yu2016decomposing,Jianjie2017A}.

\item {\it{Comparison of objective functions}}: {\it{1) Weights}}: Different from the reference points, the weights (or weight vectors) show the relative importance of the objectives from the perspective of the DM \cite{5585741}, which are commonly embedded into the aggregation functions such as the weighted sum approach \cite{miettinen2012nonlinear}, Tchebycheff approach \cite{miettinen2012nonlinear}, and their modifications \cite{7572016,Branke2005}. {\it{2) Preference relation}}: Slightly different from the weight vectors, the preference relation is only a contrast relation to reflect which objective is more, equal or less important than another objective \cite{zitzler2008spam, lopez2009study, saaty2008decision} such as the analytic hierarchy process \cite{saaty2008decision}.  {\it{3) Reference vectors}}: They indicate the evolution directions of the sub-problems with respect to different objectives and usually originate from the ideal point (or the original point) to the reference points \cite{6600851}. Solutions associated with them will be guided towards the intersections between the reference vectors and $\mathcal{PF}$ such as the penalty-based intersection approach \cite{4358754,zheng2017decomposition} and angle penalized distance \cite{7386636}. {\it{4) Trade-offs}}: Trade-offs can provide the search directions on the manifold of the $\mathcal{PF}$, but the DM needs to set the amount of improvement on some or all objectives except for the reference objective which suffers from a unit of sacrifice \cite{miettinen2012nonlinear}.  {\it{5) Classification of objectives}}: The main idea is to separate different classes of objectives according to the aspirations of the objectives \cite{Miettinen2008}. For example, one class of objectives ($f_1$ and $f_2$) should be increased to meet the specified upper bounds.  

\item {\it{Comparison of solutions}}: {\it{1) Pairwise comparison of solutions}}: The key motivation is to judge the solutions in the population whether they are preferred or incomparable with each other, such as the outranking relation \cite{2002Axiomatization,siskos1986use} and preference ranking  \cite{waegeman2011era}.  {\it{2) Classification of solutions}}: The core idea is to separate the solutions into multiple classes where the solutions in the same class are incompatible with each other \cite{5585982}.  {\it{3) Selecting the most preferred solution}}: It intends to choose to best solution among the given solution set based on the qualitative preference information from the DM \cite{2008Ordinal}.
\end{itemize}

{\it{{Implicit preference:}}} In the multi-objective optimization, there are some situations where DMs cannot explicitly articulate their preferences, due to the lack of \textit{a priori} knowledge of the highly complex optimization problems, the intractable and laborious interactive process, or few available observations which result from the resource-intensive and time-consuming computation \cite{8790106}. Consequently, the implicit preference, especially the knee point \cite{GuoSurvey2021}, has attract much attention and been adopted as the preference in the multi-objective optimization \cite{8638825}. The knee points are located in the regions where bulge out the most on the $\mathcal{PF}$ \cite{2011Understanding} and they need a large compromise (or deterioration) in at least one objective to get small improvement on other objectives \cite{1999On}, as shown in Fig. \ref{plot-knees}. Interested readers are referred to the comprehensive review \cite{GuoSurvey2021} for more details about the knee points.

\section{From Preferences to Fairness}\label{SectionofP2F}
In this section, we will elaborate the relationship between preferences and fairness. Firstly, an introduction of fairness in ML is presented, including the sources of the bias, notions of the fairness, and the approaches to ensure fairness in ML. Then, the differences between the fairness and preferences are discussed. 

\begin{figure}[tb]
\centering
  \includegraphics[width=.65\linewidth]{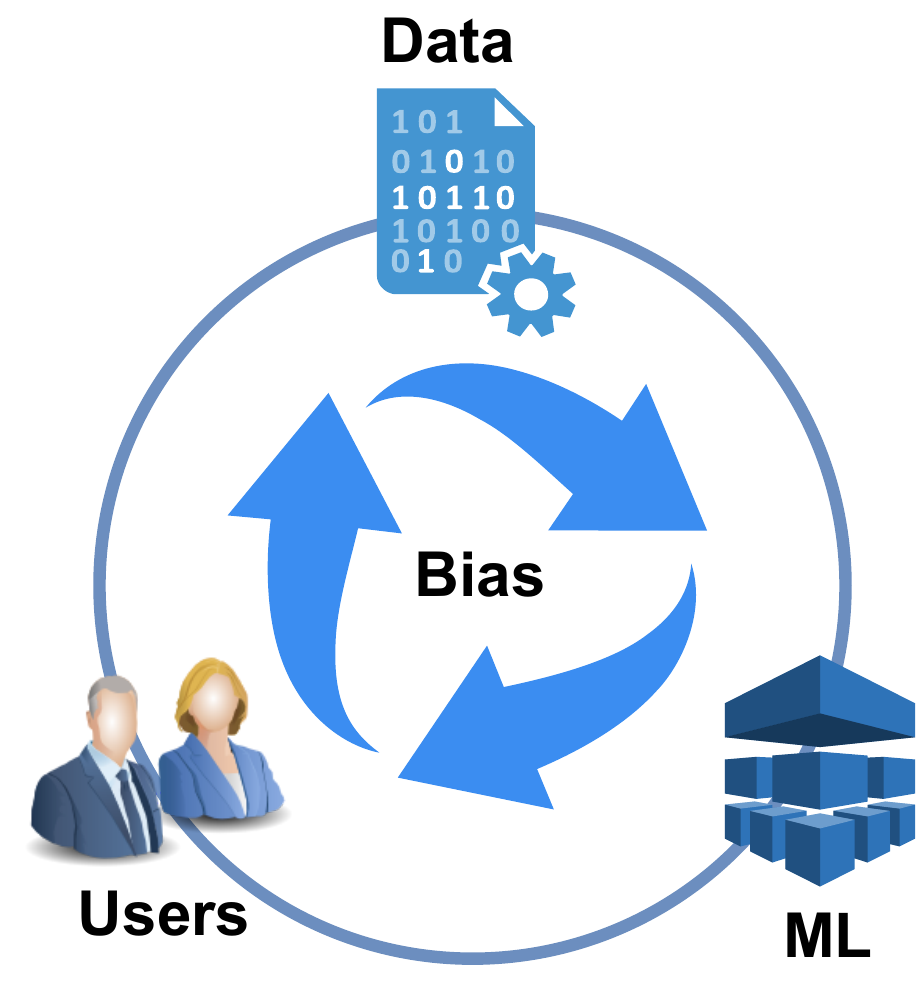}\\
  \caption{An example of the feedback loop of bias among data, ML models and users.}\label{plot-loopbias}
\end{figure}

\subsection{Fairness in machine learning}
\subsubsection{Bias in ML}
In Fig. \ref{plot-loopbias}, the feedback loop of bias (or discrimination) \cite{mehrabi2021survey} has shown that the biases are not only in data, algorithms or ML models but also in user experiences \cite{caton2020fairness,quy2021survey,olteanu2019social}.  

\begin{itemize}
\item {\it{Bias from data:}} When the training data contains biases, the ML models trained on such data will be impacted and may produce biased predictions or unfair outcomes \cite{zhang2018mitigating}. The bias in data may come from, among others, human prejudice, stereotypes based on sensitive attributes \cite{suresh2019framework}, ignorance of specific features \cite{clarke2005phantom}, or non-representative samples \cite{suresh2019framework,buolamwini2018gender}.  

\item {\it{Bias from model:}} The design of ML models may produce algorithmic bias \cite{baeza2018bias}, such as the improper parameter settings of the designed optimization methods and regularization, improper regression models, or statistically biased estimators \cite{danks2017algorithmic}. Besides, the bias may also appear in the presentation of the results, like making the popular or top-ranked items easily exposed on a website \cite{baeza2018bias,ciampaglia2018algorithmic}. In addition, inappropriate evaluation of the model may also produce unfairness or bias \cite{suresh2019framework}, such as disproportional benchmarks in facial recognition \cite{buolamwini2018gender}.

\item {\it{Bias from users:}} The bias mainly comes from the data generation or user interactions, and any inherent biases in users may lead to biased data which then further make the ML model unfair \cite{suresh2019framework}. For example, non-representative data may comes from different platforms having different user demographics \cite{olteanu2019social}. 
In addition, decisions of a user can be influenced by others' decision making (or the average group statistical results) \cite{baeza2018bias}, and users behaviors also vary from different platforms \cite{olteanu2019social}. Notably, the bias may also come from the difference in the contents generated by users in different genders, ages, regions, or with different languages \cite{nguyen2013old}. 
\end{itemize}

\subsubsection{Notions of fairness in machine learning}
According to the surveys \cite{olteanu2019social,mehrabi2021survey,caton2020fairness}, the fairness metrics can be roughly classified into group-based fairness, subgroup-based fairness, and individual-based fairness. 
\begin{itemize}
\item {\it{Group-based fairness}} aims to make different groups equally treated. Several metrics of group-based fairness are introduced as follows. {\it{Demographic parity}} \cite{dwork2012fairness} denotes that each group has equal probability to be classified with the positive label \cite{kamishima2012fairness,zemel2013learning}. {\it{Equal Opportunity}} \cite{hardt2016equality} indicates that the true positive rates are the same across different groups. {\it{Equalized Odds}}  \cite{hardt2016equality,berk2021fairness} means that the protected and unprotected groups have the same true positive rate and false positive rate. {\it{Treatment equality}} \cite{berk2021fairness} means that the ratio of false negative prediction to false positive prediction is same for both protected and unprotected groups. {\it{Test fairness (calibration)}} \cite{chouldechova2017fair,verma2018fairness} aims to ensure people have the same probability to be assigned to a positive class when they are from different groups but with the same predicted probability score. {\it{Conditional statistical parity}} \cite{corbett2017algorithmic} means that different groups of people have equal probability of being assigned to a positive class given a set of legitimate factors.

\item {\it{Subgroup-based fairness:}} imposes a statistical constraint across a large number of {\it{subgroups}} defined by the protected attributes. For example, {\it{Subgroup fairness}} picks a statistical fairness constraint (e.g., equalizing false positive rates across protected groups) and then asks the constraint whether it is able to hold over an exponentially or infinitely large collection of subgroups \cite{kearns2018preventing,kearns2019empirical}.

\item {\it{Individual-based fairness}} attempts to make similar individuals have similar predictions. {\it{Fairness through unawareness}} \cite{grgic2016case} aims to ensure that any protected attributes or features are not explicitly used in the decision-making process \cite{kusner2017counterfactual}. {\it{Fairness through awareness}} \cite{dwork2012fairness} indicates that similar individuals should have similar predictions.
{\it{Counterfactual fairness}} \cite{kusner2017counterfactual} ensures that a decision is fair towards an individual if the decision in actual world coincides with that in counterfactual world where the individual belongs to a different demographic group.
\end{itemize}

\subsubsection{Mitigating unfairness in machine learning}
The fairness-aware methods in ML are commonly classified into {\it{pre-processing}}, {\it{in-processing}}, and {\it{post-processing}} approaches. 

\begin{itemize}
\item {\it{Pre-processing approaches}} remove the bias of the data or modify the features spaces unrelated with sensitive attributes before the ML model is trained. The commonest way is to learn a new representative of the data \cite{zemel2013learning,calmon2017optimized,creager2019flexibly}, which preserves as much information as possible and is insensitive to the protected attributes. In addition, some studies try to use selective weights on the tuples (attributes, label) of the training data to mitigate unfairness \cite{kamiran2012data}, or repair the feature values of data in order to increase group fairness while not breaking their ranking orders within groups \cite{feldman2015certifying,bellamy2019ai}. Besides, the re-balance of the data \cite{chawla2002smote} is essential in building fair models, which aims to to under- or over-sample observations from certain sensitive feature groups to balance the data distribution \cite{dixon2018measuring}. 

\item {\it{In-processing approaches}} change the learning algorithms by enforcing fairness constraints during the model training \cite{d2017conscientious}. Commonly, the objective of maximizing the prediction accuracy is optimized subject to the fairness constraints, or vice versa \cite{berk2017convex,bellamy2019ai,cotter2019two}. In addition, adding regulation terms of fairness to the optimization objective is another strategy to mitigate the unfairness \cite{kamishima2012fairness}. Recently, multi-objective evolutionary learning strategy is adopted to handle such multi- or many-objective optimization problem with respect to the model accuracy and different fairness metrics \cite{fairevo2021,liu2022accuracy}. Adversarial learning is also used to mitigate unfairness by maximizing the prediction ability and minimizing the adversary's ability in predicting protected attributes \cite{zhang2018mitigating}.

\item  {\it{Post-processing approaches}} apply transformations to the pre-trained model so as to make the outcomes fair on the condition of given fairness metrics. For example, calibrating decision threshold is a common way to mitigate the unfairness \cite{hardt2016equality,fish2016confidence}. Moreover, Pleiss et al. \cite{pleiss2017fairness} came up with a relaxed equalized odds with calibration to find probabilities for the balance between achieving parity and preserving calibrated probability estimates. Kamiran et al.\cite{kamiran2012decision} proposed a reject option classification to find the balance between 
providing positive labels to deprived groups and negative labels to favored groups under uncertainty.
\end{itemize}

\subsection{Differences between preference and fairness}\label{SecDpf}
\begin{figure}[t]
\centering
  \includegraphics[width=.9\linewidth]{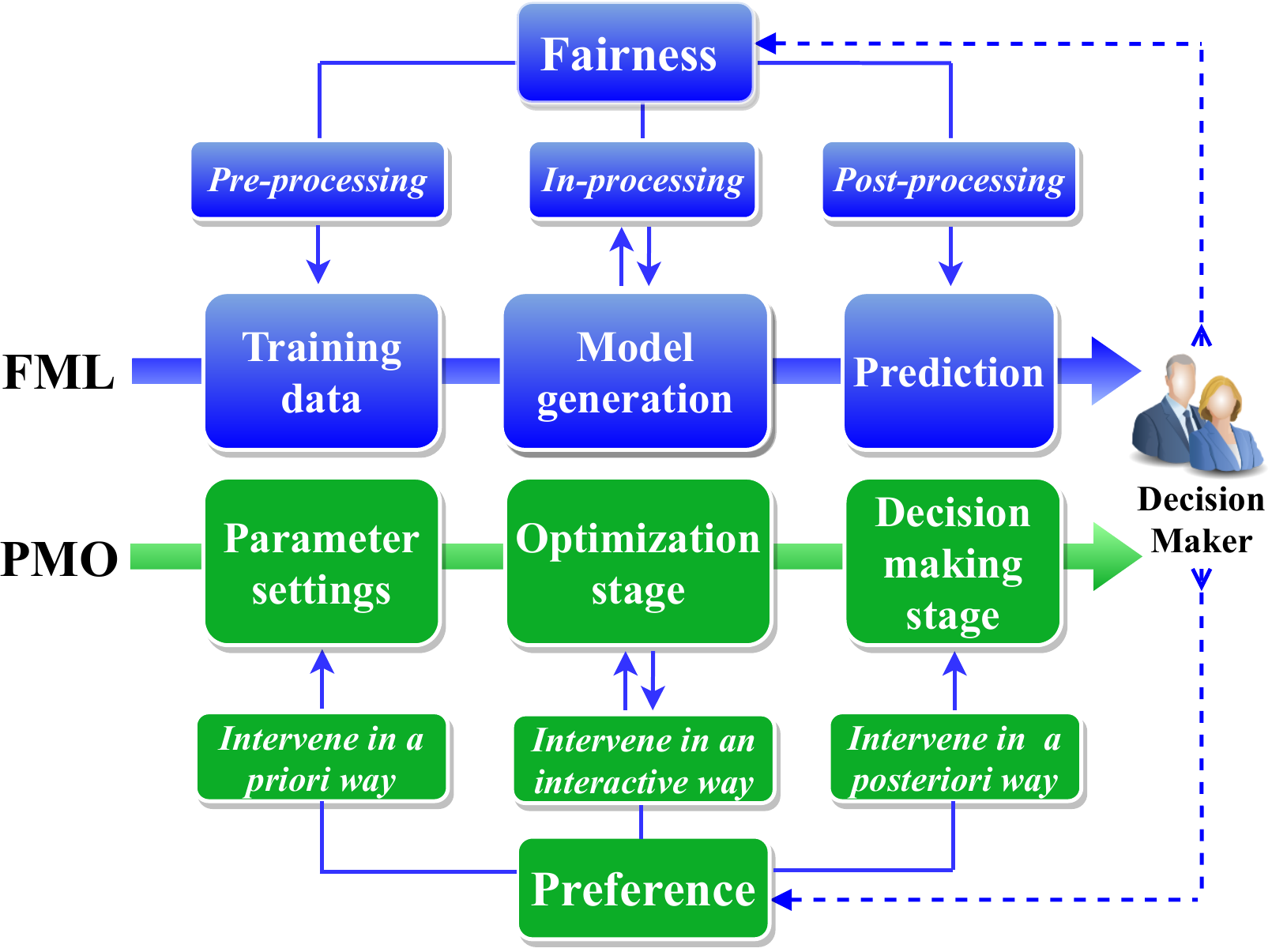}\\ 
  \caption{Comparison between the PMO and FML, in both of which the user can interact before, during or after the optimization or learning process. }\label{plot-P2F}
\end{figure}

As shown in Fig. \ref{plot-P2F}, we can see that PMO and FML show many similarities and differences.  

\subsubsection{Similarities}
The fairness metrics described in the previous subsection are kind of preferences from DMs. Different backgrounds of culture or other social factors make DMs to have different perspectives in looking at the fairness \cite{1993TradeOffs}. As a result, no single definition of fairness can satisfy all DMs, and different fairness metrics may even conflict with each other \cite{chouldechova2017fair,zhang2020joint}. For example, there is a trade-off between the individual fairness and group fairness \cite{chouldechova2017fair}. In fairness-aware federated learning, the collaboration fairness ensures that the clients with more contributions will get higher rewards \cite{shi2021survey}. Therefore, different definitions or descriptions of fairness potentially represent the expectations of different groups of people. 

The unfairness or inequality sometimes comes from the individuals' preference \cite{avin2015homophily,kim2019preference}. For example, the social inequality and discrimination in social networks arises from predefined communities and people's preferences in building new connections \cite{avin2015homophily}. An underlying assumption of FML is that the preference of the DMs is often known beforehand. Although it does not always hold in practice \cite{finocchiaro2021bridging}, preferences can be considered in the preference-driven FML to improve model accuracy without violating the preferences \cite{kim2019preference,do2021online,zafar2017parity}. Taking complex online advertising as example, we can see that the targeted advertising is in line with the preference of the DM,  which is often unknown in advance and will be estimated in practice \cite{kim2019preference,liu2019personalized}.

\subsubsection{Differences}
In a broad sense, fairness is the absence of any bias towards a DM or groups based on their sensitive attributes in the context of decision-making \cite{chouldechova2017fair,saxena2019fairness}. However, Kim et al. \cite{kim2019preference} show that less-preferred outcomes will be found under the constraint of individual fairness when DMs have their preferences over the outcomes. In other words, preferences and fairness are not necessarily consistent with each other. In a narrow sense, preferences are commonly articulated as threshold to determine the tolerance of unfairness  \cite{fish2016confidence,kamiran2012decision,menon2017cost,valera2018enhancing}, or the trade-off preference between the fairness and model accuracy \cite{fairevo2021,liu2022accuracy}.  
That is, preferences in ML are meant for assisting in finding acceptable fair outcomes for DMs. 

Since the outcomes finally serve the DMs, a number of studies have been proposed to satisfy the DMs' preferences over the outcomes of the decisions \cite{kim2019preference}. For example, Zafar et al. \cite{zafar2017parity} proposed two preference-based notions of fairness to allow any group of DMs to collectively prefer the outcomes without considering the fairness as compared to the other groups. Their results have shown that the preference-based fairness often achieves higher accuracy than other parity-based fairness notions. Valera et al. \cite{valera2018enhancing} have revealed that the accuracy and fairness can be significantly improved by estimating the thresholds used by DMs to balance the exploitation (selecting accurate and fair decisions) and exploration (selecting decisions satisfying DMs' preference).

\section{Fairness in multi-objective optimization}\label{SectionofFMO}
In this section, we will introduce notions of fairness in multi-objective optimization from various perspectives. Then, the differences between the fairness in multi-objective optimization and ML will be discussed.

\subsection{Fairness in traditional multi-objective optimization}\label{SectionFTMO}
\subsubsection{Fair sampling in multi-objective optimization}
In 2004, Laumanns et al. \cite{laumanns2004running} came up with theoretical running time analysis on multi-objective evolutionary algorithms from the aspect of the fairness of sampling. By {\it{fair sampling}}, it is meant that the number of offspring generated by all individuals are equal at the end of optimization. Their experimental results show that {\it{fair sampling}} is able to accelerate the optimization. Following the above idea, Friedrich et al. \cite{friedrich2008runtime,friedrich2011illustration} took further investigation on how the fairness influences the runtime behavior by considering the {\it{fair sampling}} in decision space or objective space. Their experimental results have demonstrated that the influence from different {\it{fair sampling}} mechanisms varies significantly with the runtime of the optimization process on different optimization functions. Qian et al. \cite{qian2016selection} also investigated the effectiveness of {\it{mixing fair selection mechanisms}} (w.r.t. the decision space or objective space) to multi-objective optimization, where the fairness is to ensure the balance of number of offspring of all individuals in current population. Their results indicate that the advantages of one fair selection mechanism can compensate for the disadvantages of another.

\subsubsection{Fairness in the optimization of dynamic prices for electric vehicle charging}
Dynamic pricing scheme is one of popular methods for the reduction of the cost of electric vehicle charging. The optimization of dynamic prices can be described as an MOP  by optimizing the prices to maximize the profit of charging station operator, minimize the number of declining customers either when the price is too high or when the stations are occupied \cite{limmer2017multi}. However, varying prices may cause unfairness that people experience different prices for similar products \cite{2004TheUnfairPrice}. Hence, Limmer et al. \cite{limmer2018optimization} consider two  notions of fairness and transfer them into two objectives in the optimization of dynamic prices. One is that customers should receive similar prices when they arrive in the same charging time window requesting similar orders. The other is that customers should get similar price offers when they have similar requests arriving in different but nearby time windows. Their experimental results show that the fairness-aware optimization method is able to increase the fairness and keep the expected profit, but may increase the number of declining customers. To some extent, the maximization of fairness does not impact the profit. 
However, more investigations are essential on the fairness for different charging strategies \cite{lunz2012influence}.

\subsubsection{Fairness in optimization of next release problem}
The next release problem (NRP) \cite{bagnall2001next} is to find the balance between the cost and benefits for the next release of software under a number of requirements from customers. It often happens that requests from a customer are inconsistent with those from another. Therefore, selecting appropriate requirements for the release of next version of software is a critical decision to make for software companies \cite{dong2022multi}. 

Finkelstein et al. \cite{finkelstein2009search} investigate the trade-off in different notions of fairness in the optimization of NRP, in order to find the fair allocation of requirements to the clients. To this end, they have explored the {\it{fairness in requirement assignment}}. Three factors to balance the fulfilled requirements between the clients are taken into account, i.e., the number, value and cost of the requirements. This study has proposed three models of fairness:  {\it{fairness on absolute number of the fulfilled requirements}}, {\it{fairness on absolute value of fulfilled requirements}}, and  {\it{fairness on the percentage of value and cost of fulfilled requirements}}. Their experimental results show that there is a strong trade-off between the number of requirements and fairness to the clients. Their study can help DMs provide the clients with fair allocation of requirements based on various fair metrics. During the decision-making stage, Ghasemi et al. \cite{ghasemi2021multi} suggested to take different quality indicators like {\it{fairness in requirement assignment}} to compare the obtained Pareto non-dominated solutions of NRP \cite{finkelstein2009search}, which plays an important role in requirement analysis or decision making. 

\subsubsection{Fairness in preference learning}
Preferences are commonly to show a partial ordering over outcomes, and the aggregation of  the pair-wise preference to a final decision is crucial. However, different DMs may have varying voting rules, while ensuring fairness of the aggregation of the pair-wise preference attracts much attention. Rossi et al. \cite{rossi2005aggregating} define the fairness of preference aggregation from four aspects. {\it{Freeness}} is to ensure no restriction on the outcomes. {\it{Independence to irrelevant alternatives}} is to ensure the relation between two outcomes (or alternatives) depending only on their preference relation given by the DMs (and not on their preferences over other outcomes). {\it{Monotonicity}} is to ensure that moving up the position of an outcome in one preference ordering does not lead to the moving-down of the position of the outcome in other preference ordering. {\it{Non-dictatorial}} is to ensure that there is no DM such that, no matter what other DMs give their preferences, the DM is never contradicted in the outcome. They have studied a number of preference aggregations and demonstrated that under certain conditions on the partial orders, no preference aggregation system can be fair when there are at least two DMs and three outcomes to order.
The conclusion has also been proved in \cite{zhang2005fairness} which indicates the unfairness comes from the dependence to irrelevant outcomes.

\subsubsection{Fairness in resource allocation}
Resource allocation problems (RAPs) \cite{luss1999equitable}  are ubiquitous in many areas, such as task scheduling \cite{bowerman1995multi,purushothaman2021evolutionary}, emergency service allocation \cite{marsh1994equity,jagtenberg2021introducing} and Cloud service allocation \cite{JIANG2016239,6211077}.  In recent decades, equitable allocation of resources that ensures equitable distribution of resources to users \cite{FairCAS2007,fairRA2021,ghodsi2011dominant}, has become a popular research topic. 

In handling the network resource allocation problem, Michael et al. \cite{2002Linear} have theoretically studied the linear optimization of RAPs with respect to multiple equitable criteria. They have proposed the notion of {\it{equitable efficiency}} to identify the preferred solutions having better distribution of the outcome values (i.e., objectives of an MOP) from a large number of Pareto optimal solutions. The principle behind the {\it{equitable efficiency}} is that any small amount of transfer from an
objective to any other relatively worse-off objective produces a more preferred solution in the objective space. For example, a solution $B$=$(2,2,2)$ with better distribution on objective values is generated from solution $A$=$(4,2,0)$ after two units transfer from $f_1$ to $f_3$. 
Following the above principle, {\it{equitable rational preference relations}} \cite{kostreva2004equitable} is developed to refine the selection of Pareto optimal solutions, which is also applied in network resource allocation \cite{ogryczak2014fair,koppen2010comparison,le2005rate}. 

Similarly, there are also increasing concerns about fairness on facility location problems (FLPs) that aim to find optimal locations of facilities to satisfy the demand of customers.   
In \cite{blanco2022fairness}, the notion of fairness is to ensure that the maximization of the covered demand of one facility does not negatively impact the demand coverage of others. In ambulance location problem \cite{jagtenberg2021introducing,jagtenberg2020fairness}, the fairness is denoted by the maximization of {\it{Bernoulli-Nash social welfare}} \cite{burlacu2020solving}, i.e., maximizing the joint probability that everyone receives ambulance on time.  

In addition, notions of fairness are also considered in the water resource allocation problem, since the overuse or poor management of water resources may lead to droughts or water crisis \cite{habibi2016multi}. In handling intermittent operation of the water distribution networks problem, Solgi et al. \cite{solgi2020multi} formulate the fairness into the uniform distribution of water supply to the water consumption nodes of the network. In \cite{fu2018water,fu2021comparison}, the equity principle denotes the fair use of water resources, i.e., the users having higher levels of water use are allocated with more water resources. In addition, Tang et al. \cite{tang2021new} have proposed a multi-objective optimization model to deal with the water resource allocation problem. The designed objectives are to minimize water shortage risks and ensure the fairness of regional water distribution, where the fairness is described by {\it{Gini coefficient}} \cite{cullis2007applying} across agricultural, industrial and domestic sections.

\subsection{Fairness in data-driven optimization}
Data-driven optimization aims to solve the optimization problems where the evaluation of objective (or fitness) or constraints are based on data collected from physical experiments, numerical simulations, or everyday life \cite{8456559,jin2021data,li2022evolutionary}. Recent years have seen a few single-objective cases of fairness-aware data-driven optimization, while no fairness-aware multi-objective optimization has been reported.

\subsubsection{Data-driven turbulence modelling}
In \cite{jiang2021interpretable}, Jiang et al. take into account the fairness in the data-driven turbulence modelling for engineering simulations, where the fairness is to ensure no prejudice in predicted outcomes across multiple classes of flows. To this end, this study has proposed {\it{a clustering-based unbiased data sampling}} and {\it{a fair cost function design}} since the unfairness may come from the training data or cost functions. The clustering-based unbiased data sampling is based on an unsupervised classifier without \textit{a priori} knowledge, which is able to re-balance the sample diversity between different subregions. In other words, the fair sampling technique is beneficial to the development of the data-driven turbulence model across different flow patterns. The latter ensures the applicability of the hyperparameter settings of the training model across similar cases. In addition, it guarantees fair contributions of all stress anisotropy components of turbulence model to the overall training error. Their results show that the training model has good generalization across two- and three-dimensional flows.

\subsubsection{Data-driven wireless resource allocation}
In \cite{sun2022learning}, Sun et al. have proposed a new memory-based continual learning formulation for wireless system design, which enables data-driven approaches to continuously learn and optimize wireless resource allocation in a dynamic environment. In this study, a  deep neural network is trained based on current data batch and a working memory of representative historical data for dynamic environment. Especially, {\it{a notion of sample-fairness criterion}} \cite{hong2014signal} is introduced to give the data samples having relatively lower performance higher probability to be selected into the memory. The assumption is that the model will perform well on data samples in dynamic environment if the model is capable of achieving good performance on challenging and under-performing samples. Their experimental results have demonstrated that the fairness-aware data-sample selection is able to contribute to the competitive performance of the data-driven approach for wireless resource allocation in dynamic environment. 

Referring to fair resource allocation, Kaur et al. \cite{kaur2022vaxequity} have developed data-driven pandemic risk assessment and optimization framework to achieve fair COVID-19 vaccine distribution among COVAX\footnote{COVAX is an initiative from WHO that supplies vaccines to the most needed countries.} countries. Their results have shown that the data-driven risk assessment model is able to facilitate the optimal distribution of limited vaccine resource in the optimization of the Jain's fairness \cite{5461911} embedded objective function, where the fairness ensures that the vulnerable countries lacking vaccine resources have a priority to get support under the consideration of combined risk across all COVAX countries. 

\subsubsection{Data-driven fairness-aware Bayesian optimization}
Bayesian optimization is a popular and effective data-driven optimization approach \cite{wang2022survey}. In \cite{perrone2021fair}, Perrone et al. have proposed a general constrained Bayesian optimization framework which is able to cater different ML models and one or multiple fairness constraints. Following the Bayesian optimization \cite{williams2006gaussian}, this study iteratively tunes hyperparameters ($\mathbf{x}$) based on the best query ($\mathbf{x}^*$), which is achieved by maximizing an acquisition function. The acquisition function is built on a posterior Gaussian process with respect to the designed objective $f(\mathbf{x}$) (i.e., accuracy) and fairness constraint $c(\mathbf{x})$. The fairness in this study is expressed as the probability of satisfying the fairness constraint $P(c(\mathbf{x})\leq \epsilon)$ \cite{gardner2014bayesian,gelbart2014bayesian}, where $\epsilon$ is an unfairness bound \cite{donini2018empirical}. Notably, given $D$ fairness constraints, then a product of the probabilities is replaced, i.e., $P=\Pi^D_{d=1}P(c_d(\mathbf{x})\leq \epsilon_d)$. Once a new query ($\mathbf{x}^*$) is found, then it will be evaluated by the designed objective and constraint functions, i.e., ($f(\mathbf{x}^*)$, $c(\mathbf{x}^*)$). After that, the pool of $\{\mathbf{x}, f(\mathbf{x}), c(\mathbf{x})\}$ is updated with the evaluated solution. After a number of iterations, the best fair hyperparameter configuration will be acquired from the pool of $\{\mathbf{x}, f(\mathbf{x}), c(\mathbf{x})\}$. Their experimental results show that the proposed fairness-aware Bayesian optimization is flexible to arbitrary fairness measures and competitive with the state-of-the-art algorithmic fairness techniques.

\begin{figure}[t]
\centering
  \includegraphics[width=.85\linewidth]{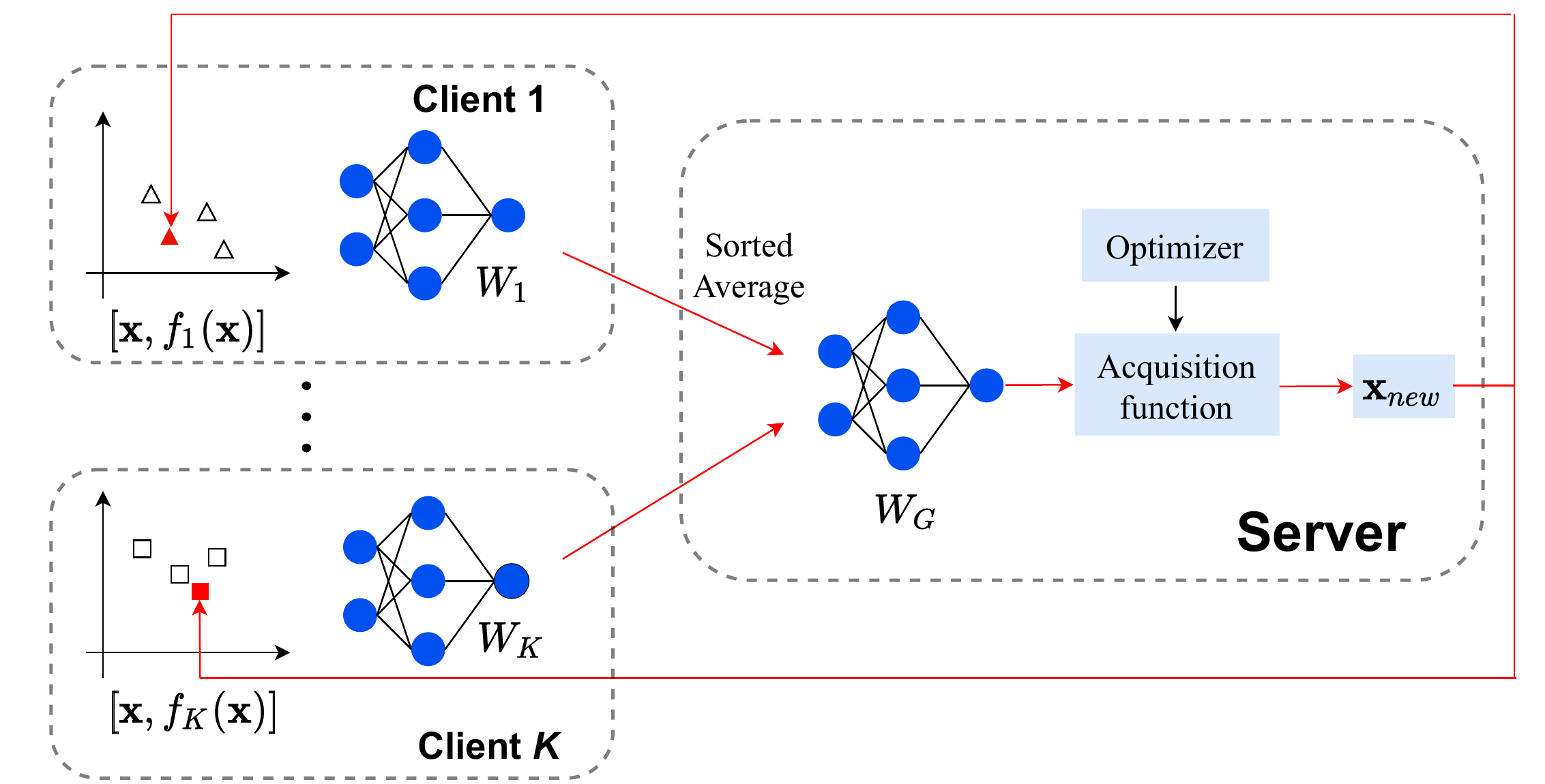}\\
  \caption{A generic framework of federated optimization. The global model $W_G$ is obtained by conducting sorted average on the local models ($W_1,\ldots W_K$). The acquisition function is formulated on the basis of the global and local models and optimized by any optimizer to get a new candidate solution ${\bf{x}}_{new}$ for the update of local models.
}\label{plot-FOP}
\end{figure}

\subsection{Fairness in federated optimization}\label{Sec-FFO}
Federated optimization is firstly proposed in  \cite{xu2021federatedEA}, where the data-driven optimization collects data in a distributed way when the data is subject to privacy restrictions.
In Fig. \ref{plot-FOP}, an example of the federated optimization is given. Undoubtedly, federated optimization is a special case of data-driven optimization \cite{jin2021data,8456559,li2022evolutionary}. The main difference lies to the global model ($W_G$) obtained by sorted averaging on local models ($W_1,\ldots W_K$) to approximate the fitness functions (commonly refer to the objectives) \cite{jin2021data}. Another difference is the training data distributed on the local clients and cannot be directly manipulated in an ensemble way. The fundamental difference between the federated optimization and federated learning is that federated optimization aims to assist the optimization process in finding the global optima or $\mathcal{PF}$ of the corresponding single optimization problem  \cite{xu2021federatedEA} or MOP \cite{xu2021federatedMOEA}, respectively, while federated learning is to train highly accurate global model \cite{zhu2021real,2021FromFL}.

From Fig. \ref{plot-FOP}, we can see that federated optimization shares many similarities with federated learning such as the similar federated framework with distributed data, and model aggregation. Unfortunately, there are few studies on the fairness-aware federated optimization, although the fairness-aware federated learning has been fast developed in recent years \cite{shi2021survey}. Especially, we can find out that the fairness issues may also exist in federated optimization. For example, the accuracy of the local models of federated optimization varies because a local model may not be updated when the new candidate solution is located in the infeasible region of optimization problem \cite{xu2021federatedEA,xu2021federatedMOEA}. As a result, inaccurate models may mislead the optimization process. Accordingly, the {\it{performance fairness}} \cite{yue2021gifair,li2019fair} should be taken into account that uniform accuracy distribution across the local models should be ensured. Another example is that the clients are not willing to share their data without incentive rewards. As a result, the {\it{collaboration fairness}} \cite{lyu2020collaborative} can be borrowed from fairness-aware federated learning that the participants with higher contribution will receive higher rewards \cite{shi2021survey}. In federated optimization, the rewards may be described as the preference for some objectives (or local models). In addition, the training data with sensitive attributes may make the accuracy of the local models various, leading to biased optimization on different objectives. Therefore, the biased optimization may result in the solutions located in a small region of the $\mathcal{PF}$ that only satisfies certain groups of people. In other words, biased optimization may harm the benefits of specific groups of people who are not interested in the obtained ROIs. A fair example is given in Fig. \ref{plot-prefer} where the ROIs consist of the solutions under diverse preferences from different groups of people. Therefore, the notion of {\it{model fairness}} \cite{chouldechova2018frontiers,mehrabi2021survey,barocas-hardt-narayanan} from FML can be adopted to handle the fairness-aware federated optimization problem to ensure that no discrimination of the model against specific groups of people will happen.
 
\subsection{Differences between fairness in optimization and machine learning}\label{DiffFinOpML}
ML is to learn or abstract a model and its corresponding hyperparameter configurations based on the given data set, while the fairness in ML mainly acts on the relationship between the learned model and abstracted hyperparameter configurations \cite{ntoutsi2020bias}. Most fairness-aware machine learning methods aim to ensure that the model performance is insensitive to the protected or sensitive attributes  \cite{chouldechova2018frontiers,mehrabi2021survey,zehlike2021fairness}. By contrast, the fairness considered in optimization is to ensure the equity among the individuals (or groups) and there is no concern with the sensitive attributes. Therefore, fairness in optimization may lie in all aspects influencing the optimization process, such as the {\it{fair-sampling}} \cite{friedrich2008runtime,friedrich2011illustration}, the design of fairness-aware objective \cite{limmer2018optimization}, fairness between the preferences from DMs \cite{bagnall2001next,rossi2005aggregating}, and {\it{fair-outcomes}}  \cite{purushothaman2021evolutionary,jagtenberg2021introducing,sun2022learning}. 

In essence, all machine learning problems are optimization problems \cite{2006Multi} where the optimization objectives are typically designed with respect to the learning model, such as the maximization of model performance (e.g., accuracy), maximization of diversity \cite{gu2015ensemble}, minimization of model complexity \cite{2008Pareto}, maximization of communication efficiency in federated learning \cite{zhu2021real,2021FromFL}, maximization of robustness \cite{liu2021NC,liu2022TEVC}, among many others \cite{2006Multi,wissam2015cyb}. Specifically, fairness metrics are regarded as extra objectives in the optimization to ensure the fairness of the learning model \cite{fairevo2021,zhang2020joint,valdivia2021fair}. Different from conventional learning methods to solve one of the objectives or a scalar function of multiple objectives, the Pareto-based multi-objective machine learning \cite{2008Pareto,EvoML-Survey,ZHAN202242} simultaneously optimizing multiple objectives has received increased interest recently. The reason is that we can gain a deeper insight into ML problems by analysing the $\mathcal{PS}$ and $\mathcal{PF}$. Therefore, from the above perspective, fairness-aware machine learning can also be seen as a category of fairness-aware optimization. 

\section{Challenges and opportunities}\label{SectionofCO}
According to the above discussions, there are several promising lines of research embodied with a number of challenges and opportunities. 

\subsection{Fairness-aware multi-objective optimization}
According to the discussions in Section \ref{SectionofFMO}, the following lines of research are promising and of great importance to handle challenges in FMO.

\begin{itemize}
 \item As discussed in Section \ref{SecDpf}, different definitions of fairness may satisfy different groups of people due to different backgrounds of culture or other social factors  \cite{1993TradeOffs}. Accordingly, one promising line of research is how to satisfy the preference of the DMs to design proper fairness metrics to different applications \cite{limmer2018optimization}. Another common issue to be solved is whether we should take fairness as objectives, constraints, or penalize objective with fairness requirements.
 
\item It has been demonstrated that there is a trade-off between the model performance and fairness \cite{fairevo2021,zhang2020joint,valdivia2021fair} and it is challenging to mitigate unfairness among all clients or groups of people. Hence, more investigations may be essential on integrating preferences from the users into the fairness-aware optimization, in order to alleviate the issue of unfairness \cite{ruchte2021scalable,ustun2019fairness}
 
\item It has been found theoretically that fair sampling may affect the runtime of the optimization process \cite{friedrich2008runtime,friedrich2011illustration,qian2016selection}. However, more investigations are required on the fair performance indicators \cite{8314455} and fair comparison between the outcomes \cite{rossi2005aggregating}.

\item There is a specific kind of expensive optimization problems where different objectives may have different evaluation time or number of evaluations \cite{wang2021transfer,wang2021transfer2}. Hence, how to ensure the fairness between the optimality of the data-driven objectives is essential. On the contrary, a question arises that all objectives are equally important or not during the optimization  \cite{7088618}. 


\end{itemize}

\subsection{Fairness-aware data-driven optimization}
Generally speaking, there are two lines of research in the fairness-aware data-driven optimization. One is to build data-driven models or surrogates to describe the fairness when the fairness cannot be explicitly articulated. The other is to find a fair model for handling the ML tasks when the fairness can be explicitly expressed. Accordingly, we can list a number of challenges of the data-driven optimization according to  \cite{ntoutsi2020bias}. First, it is challenging to accurately describe the fairness with surrogate models when the fairness cannot be explicitly articulated. Second, there are no clear conclusive results on the best choice of bias-mitigating method on different category of tasks in different intervention processes. Third, the fairness-aware data-driven optimization or ML is not universally valid, so that their applicable situations need to be explored like training across heterogeneous data sources, decision making online and dynamic optimization for uncertain or open environment. Fourth, the notions of fairness may be variant in different ML frameworks, such as the reinforcement learning \cite{jabbari2017fairness}, federated learning \cite{shi2021survey}, supervised learning \cite{hardt2016equality}, unsupervised learning \cite{9541160}, and deep learning \cite{du2020fairness,fitzsimons2019general}. Fifth, except the data protection law, other policies or standards to data equality or selection are essential for the development of fairness-aware data-driven optimization and ML. 

On the other hand, if we focus on their difference, more investigations may be needed on how the issues of unfairness impact the data-driven optimization process. For example, the biased  sampling process may generate different distributions of training samples, which may impact the accuracy of surrogate models \cite{friedrich2008runtime,friedrich2011illustration}. Another interesting topic is to learn preferences from DMs to set the appropriate level of fairness. For example, some researches take the probability of satisfying the fairness constraint $P(c(\mathbf{x})\leq \epsilon)$ \cite{gardner2014bayesian,gelbart2014bayesian} as the consideration of fairness, where $\epsilon$ is an unfairness bound \cite{donini2018empirical} and $\epsilon$ can be learnt from the DMs' preference. In addition, fairness between the preferences from DMs \cite{bagnall2001next,rossi2005aggregating} is essential for fair decision making when multiple DMs are involved. 

\subsection{Fairness-aware federated optimization}
Fairness-aware federated optimization will be an emerging research topic since it shares many similarities with the fairness-aware federated learning \cite{shi2021survey}. The {\it{performance fairness}} \cite{yue2021gifair,li2019fair}, {\it{collaboration fairness}} \cite{lyu2020collaborative}
and {\it{model fairness}} \cite{chouldechova2018frontiers,mehrabi2021survey,barocas-hardt-narayanan} in fairness-aware federated learning may also appear in fairness-aware federated optimization as introduced in Section \ref{Sec-FFO}.
In addition, {\it{fair sampling}} may be another line of research. When the distribution of the training data is biased, the surrogate models trained on the data sets in different local clients may lead to biased optimizations on different objectives after the average aggregation. In \cite{jiang2021interpretable}, it has been proved that the accuracy of the surrogate model is sensitive to the data distribution. 
Besides, it remains to be investigated that the average aggregation of local models is fair or not, since different clients may have varying voting rules \cite{rossi2005aggregating}. Furthermore, how to ensure the fairness between different local models is an issue, when the local models have totally different structures or different types of optimization problems are in different clients. For example, the optimization problem on a client is a single-objective optimization problem, while the problems on others are multi- or even many-objective optimization problems.    

\subsection{Fairness-aware multi-objective machine learning}
\subsubsection {{Fair machine learning through multi-objective optimization}} 
In past decades, traditional methods for FML commonly either optimize the utilities of ML but subject to the constraints of fairness \cite{MatthewNips2016,singh2018fairness,zehlike2017fa,komiyama2018nonconvex} or optimize the fairness but under the limit of lower bound of utilities
\cite{biega2018equity}. Recent, increasing attention has been paid to the optimization of both utilities and fairness through multi-objective optimization \cite{fairevo2021,valdivia2021fair,schmucker2020multi,padh2021addressing,martinez2020minimax}. The studies \cite{fairevo2021,padh2021addressing,schmucker2020multi} aim to find the balance between accuracy and multiple fairness metrics, while the study \cite{padh2021addressing} uses the relaxations of fairness notions, the research \cite{schmucker2020multi} difference in fairness \cite{donini2018empirical}. Valdivia \cite{valdivia2021fair} comes up with a multi-objective framework to ensure the fairness in decision tree models with a small sacrifice of the classification accuracy. In some domains where require high-quality service such as healthcare, Martinez et al. \cite{martinez2020minimax} argue improving minimax group risk in classification is critical, so that they formulate group fairness into an MOP where each sensitive group risk is regarded as an optimization objective. As a result, the worst-case classification errors can be reduced. 
As mentioned above, the conflicts between the fairness metrics and utilities deserves more investigations since some fairness metrics themselves may or may not conflict with each other \cite{chouldechova2017fair}. In other words, properly choosing the metrics is crucial in building non-redundant fairness-aware models. 

\subsubsection {{Preference-driven fairness-aware deep multi-objective learning}}
In \cite{zhang2020joint}, Zhang et al. propose a decision framework for joint optimization of fairness and utility, which adheres to the DM's preferences. The preferences are expressed by weights to make trade-off between the fairness and utility.
Ruchte et al. \cite{ruchte2021scalable} sequentially transform the deep multi-objective learning into a single-objective optimization by inputted preference vector in each optimization iteration, and a penalty term is designed to ensure the well-spread of the preference vectors in the objective space. As a result, the whole $\mathcal{PF}$ will be approximated. Besides, a fairness objective based on a relaxation of difference of equality of opportunity is defined in their study. Their experimental results show that 
the preference-driven fairness-aware deep multi-objective learning is competitive in getting a good trade-off between all objectives. However,  preference articulations and fairness metrics from different DMs are various, especially in handling different practical applications. When the DMs cannot explicitly give their preference in fairness-aware deep multi-objective learning, the knees \cite{GuoSurvey2021} may be an alternative to alleviate the burden of decision making for DMs \cite{8638825}. 
  
\section{Conclusion}\label{SectionofConclusion}
In comparison with the extensive concerns of unfairness of ML with respect to sensitive or protected attributes, much less attention has been paid to FMO, although FMO is also commonly seen in real life. Therefore, this paper aims to emphasize the importance of fairness in multi-objective optimization and promote the investigations on FMO. To this end, we discussed the difference between the preference and fairness in multi-objective optimization for better understanding of the notions of fairness. In addition, a number of cases from different optimization frameworks are presented and discussed, followed by a discussion of the relationship between the fairness in FMO and FML. 

A number of promising lines of research have been suggested with respect to FMO and FML. In FMO, more investigations on how to ensure fair optimization are suggested, such as fair-sampling, constructing fairness-oriented objectives, and fair-outcomes. In addition, another potential topic is to balance the preference and fairness, since preference may alleviate the issue of unfairness \cite{ruchte2021scalable,ustun2019fairness}. Especially, the fairness-aware federated multi-objective optimization deserves much attention, since the notions of unfairness in federated ML may also exist in federated multi-objective optimization. For FML, two topics, including fair machine learning through multi-objective optimization and preference-driven fairness-aware deep multi-objective learning, are of particular importance. We hope that this paper will promote studies on considering fairness in multi-objective optimization problems.


\ifCLASSOPTIONcaptionsoff
  \newpage
\fi

\normalem
\bibliographystyle{IEEEtran}
\bibliography{Refer2fair}

\begin{thebibliography}{100}
\providecommand{\url}[1]{#1}
\csname url@samestyle\endcsname
\providecommand{\newblock}{\relax}
\providecommand{\bibinfo}[2]{#2}
\providecommand{\BIBentrySTDinterwordspacing}{\spaceskip=0pt\relax}
\providecommand{\BIBentryALTinterwordstretchfactor}{4}
\providecommand{\BIBentryALTinterwordspacing}{\spaceskip=\fontdimen2\font plus
\BIBentryALTinterwordstretchfactor\fontdimen3\font minus
  \fontdimen4\font\relax}
\providecommand{\BIBforeignlanguage}[2]{{%
\expandafter\ifx\csname l@#1\endcsname\relax
\typeout{** WARNING: IEEEtran.bst: No hyphenation pattern has been}%
\typeout{** loaded for the language `#1'. Using the pattern for}%
\typeout{** the default language instead.}%
\else
\language=\csname l@#1\endcsname
\fi
#2}}
\providecommand{\BIBdecl}{\relax}
\BIBdecl

\bibitem{pessach2022review}
D.~Pessach and E.~Shmueli, ``A review on fairness in machine learning,''
  \emph{ACM Computing Surveys (CSUR)}, vol.~55, no.~3, pp. 1--44, 2022.

\bibitem{9073602}
N.~Mehrabi, F.~Morstatter, N.~Peng, and A.~Galstyan, ``Debiasing community
  detection: The importance of lowly connected nodes,'' in \emph{2019 IEEE/ACM
  International Conference on Advances in Social Networks Analysis and Mining
  (ASONAM)}, New York, NY, USA, 2019, pp. 509--512.

\bibitem{WEnips2016}
T.~Bolukbasi, K.-W. Chang, J.~Zou, V.~Saligrama, and A.~Kalai, ``Man is to
  computer programmer as woman is to homemaker? debiasing word embeddings,'' in
  \emph{Proceedings of the 30th International Conference on Neural Information
  Processing Systems}, ser. NIPS'16.\hskip 1em plus 0.5em minus 0.4em\relax Red
  Hook, NY, USA: Curran Associates Inc., 2016, pp. 4356--4364.

\bibitem{WD2019}
A.~Bose and W.~Hamilton, ``Compositional fairness constraints for graph
  embeddings,'' in \emph{Proceedings of the International Conference on Machine
  Learning}, Long Beach, California, USA, 2019, pp. 715--724.

\bibitem{zafar2017fairness}
M.~B. Zafar, I.~Valera, M.~G. Rogriguez, and K.~P. Gummadi, ``Fairness
  constraints: Mechanisms for fair classification,'' in \emph{Proceedings of
  the 20th International Conference on Artificial Intelligence and
  Statistics}.\hskip 1em plus 0.5em minus 0.4em\relax PMLR, 2017, pp. 962--970.

\bibitem{berk2017convex}
R.~Berk, H.~Heidari, S.~Jabbari, M.~Joseph, M.~Kearns, J.~Morgenstern, S.~Neel,
  and A.~Roth, ``A convex framework for fair regression,'' \emph{arXiv preprint
  arXiv:1706.02409}, 2017.

\bibitem{yang2017measuring}
K.~Yang and J.~Stoyanovich, ``Measuring fairness in ranked outputs,'' in
  \emph{Proceedings of the 29th International Conference on Scientific and
  Statistical Database Management}, New York, NY, USA, 2017, pp. 1--6.

\bibitem{biega2018equity}
A.~J. Biega, K.~P. Gummadi, and G.~Weikum, ``Equity of attention: Amortizing
  individual fairness in rankings,'' in \emph{The 41st International {ACM}
  Sigir Conference on Research \& Development in Information Retrieval}.\hskip
  1em plus 0.5em minus 0.4em\relax New York, NY, USA: Association for Computing
  Machinery, 2018, pp. 405--414.

\bibitem{shi2021survey}
Y.~Shi, H.~Yu, and C.~Leung, ``A survey of fairness-aware federated learning,''
  \emph{arXiv preprint arXiv:2111.01872}, 2021.

\bibitem{angwin2016machine}
J.~Angwin, J.~Larson, S.~Mattu, and L.~Kirchner, ``Machine bias,'' in
  \emph{Ethics of Data and Analytics}.\hskip 1em plus 0.5em minus 0.4em\relax
  Auerbach Publications, 2016, pp. 254--264.

\bibitem{jabbari2017fairness}
S.~Jabbari, M.~Joseph, M.~Kearns, J.~Morgenstern, and A.~Roth, ``Fairness in
  reinforcement learning,'' in \emph{Proceedings of the 34th International
  Conference on Machine Learning}, PMLR.\hskip 1em plus 0.5em minus 0.4em\relax
  JMLR.org, 2017, pp. 1617--1626.

\bibitem{epstein2015search}
R.~Epstein and R.~E. Robertson, ``The search engine manipulation effect (seme)
  and its possible impact on the outcomes of elections,'' \emph{Proceedings of
  the National Academy of Sciences}, vol. 112, no.~33, pp. 4512--E4521, 2015.

\bibitem{braaten2011trust}
I.~Br{\aa}ten, H.~I. Str{\o}ms{\o}, and L.~Salmer{\'o}n, ``Trust and mistrust
  when students read multiple information sources about climate change,''
  \emph{Learning and Instruction}, vol.~21, no.~2, pp. 180--192, 2011.

\bibitem{babaei2021analyzing}
M.~Babaei, J.~Kulshrestha, A.~Chakraborty, E.~M. Redmiles, M.~Cha, and K.~P.
  Gummadi, ``Analyzing biases in perception of truth in news stories and their
  implications for fact checking,'' \emph{IEEE Transactions on Computational
  Social Systems}, vol.~9, no.~3, pp. 839--850, 2022.

\bibitem{ledford2019millions}
H.~Ledford, ``Millions of black people affected by racial bias in health-care
  algorithms,'' \emph{Nature}, vol. 574, no. 7780, pp. 608--610, 2019.

\bibitem{kay2015unequal}
M.~Kay, C.~Matuszek, and S.~A. Munson, ``Unequal representation and gender
  stereotypes in image search results for occupations,'' in \emph{Proceedings
  of the 33rd Annual ACM Conference on Human Factors in Computing
  Systems}.\hskip 1em plus 0.5em minus 0.4em\relax New York, NY, USA:
  Association for Computing Machinery, 2015, pp. 3819--3828.

\bibitem{hao2019making}
\BIBentryALTinterwordspacing
K.~Hao, ``Making face recognition less biased doesn't make it less scary,''
  \emph{MIT Technology Review}, 2019. [Online]. Available:
  \url{https://www.technologyreview.com/s/612846/making-face-recognition-less-biased-doesnt-makeit-less-scary/.}
\BIBentrySTDinterwordspacing

\bibitem{berk2021fairness}
R.~Berk, H.~Heidari, S.~Jabbari, M.~Kearns, and A.~Roth, ``Fairness in criminal
  justice risk assessments: The state of the art,'' \emph{Sociological Methods
  \& Research}, vol.~50, no.~1, pp. 3--44, 2021.

\bibitem{GaoSh19}
\BIBentryALTinterwordspacing
R.~Gao and C.~Shah, ``How fair can we go: Detecting the boundaries of fairness
  optimization in information retrieval,'' in \emph{Proceedings of the 2019 ACM
  SIGIR International Conference on Theory of Information Retrieval}, ser.
  ICTIR '19.\hskip 1em plus 0.5em minus 0.4em\relax New York, NY, USA:
  Association for Computing Machinery, 2019, p. 229–236. [Online]. Available:
  \url{https://doi.org/10.1145/3341981.3344215}
\BIBentrySTDinterwordspacing

\bibitem{baeza2018bias}
R.~Baeza-Yates, ``Bias on the web,'' \emph{Communications of the ACM}, vol.~61,
  no.~6, pp. 54--61, 2018.

\bibitem{speicher2018unified}
T.~Speicher, H.~Heidari, N.~Grgic-Hlaca, K.~P. Gummadi, A.~Singla, A.~Weller,
  and M.~B. Zafar, ``A unified approach to quantifying algorithmic unfairness:
  Measuring individual \& group unfairness via inequality indices,'' in
  \emph{Proceedings of the 24th ACM SIGKDD International Conference on
  Knowledge Discovery \& Data Mining}.\hskip 1em plus 0.5em minus 0.4em\relax
  New York, NY, USA: Association for Computing Machinery, 2018, pp. 2239--2248.

\bibitem{fairevo2021}
Q.~Zhang, J.~Liu, Z.~Zhang, J.~Wen, B.~Mao, and X.~Yao, ``Fairer machine
  learning through multi-objective evolutionary learning,'' in \emph{Artificial
  Neural Networks and Machine Learning -- ICANN 2021}, I.~Farka{\v{s}},
  P.~Masulli, S.~Otte, and S.~Wermter, Eds.\hskip 1em plus 0.5em minus
  0.4em\relax Cham: Springer International Publishing, 2021, pp. 111--123.

\bibitem{chouldechova2017fair}
A.~Chouldechova, ``Fair prediction with disparate impact: A study of bias in
  recidivism prediction instruments,'' \emph{Big Data}, vol.~5, no.~2, pp.
  153--163, 2017.

\bibitem{miettinen2012nonlinear}
K.~Miettinen, \emph{Nonlinear multiobjective optimization}.\hskip 1em plus
  0.5em minus 0.4em\relax Springer Science \& Business Media:
  Berlin/Heidelberg, Germany, 2012.

\bibitem{minireview2017}
H.~Wang, M.~Olhofer, and Y.~Jin, ``A mini-review on preference modeling and
  articulation in multi-objective optimization: current status and
  challenges,'' \emph{Complex $\&$ Intelligent System}, vol.~3, pp. 233--245,
  2017.

\bibitem{GuoSurvey2021}
G.~Yu, L.~Ma, Y.~Jin, W.~Du, Q.~Liu, and H.~Zhang, ``Survey on knee-oriented
  multi-objective evolutionary optimization,'' \emph{IEEE Transactions on
  Evolutionary Computation}, 2022.

\bibitem{8790106}
G.~Yu, Y.~Jin, and M.~Olhofer, ``References or preferences – rethinking
  many-objective evolutionary optimization,'' in \emph{2019 IEEE Congress on
  Evolutionary Computation (CEC)}, 2019, pp. 2410--2417.

\bibitem{adra2007comparative}
S.~F. Adra, I.~Griffin, and P.~J. Fleming, ``A comparative study of progressive
  preference articulation techniques for multiobjective optimisation,'' in
  \emph{International Conference on Evolutionary Multi-Criterion
  Optimization}.\hskip 1em plus 0.5em minus 0.4em\relax Berlin, Heidelberg:
  Springer-Verlag, 2007, pp. 908--921.

\bibitem{zhang2020joint}
Y.~Zhang, R.~Bellamy, and K.~Varshney, ``Joint optimization of ai fairness and
  utility: A human-centered approach,'' in \emph{Proceedings of the AAAI/ACM
  Conference on AI, Ethics, and Society}.\hskip 1em plus 0.5em minus
  0.4em\relax New York, NY, USA: Association for Computing Machinery, 2020, pp.
  400--406.

\bibitem{valdivia2021fair}
A.~Valdivia, J.~S{\'a}nchez-Monedero, and J.~Casillas, ``How fair can we go in
  machine learning? assessing the boundaries of accuracy and fairness,''
  \emph{International Journal of Intelligent Systems}, vol.~36, no.~4, pp.
  1619--1643, 2021.

\bibitem{ruchte2021scalable}
M.~Ruchte and J.~Grabocka, ``Scalable pareto front approximation for deep
  multi-objective learning,'' in \emph{2021 IEEE International Conference on
  Data Mining (ICDM)}.\hskip 1em plus 0.5em minus 0.4em\relax IEEE, 2021, pp.
  1306--1311.

\bibitem{chouldechova2018frontiers}
A.~Chouldechova and A.~Roth, ``The frontiers of fairness in machine learning,''
  \emph{arXiv preprint arXiv:1810.08810}, 2018.

\bibitem{mehrabi2021survey}
N.~Mehrabi, F.~Morstatter, N.~Saxena, K.~Lerman, and A.~Galstyan, ``A survey on
  bias and fairness in machine learning,'' \emph{ACM Computing Surveys (CSUR)},
  vol.~54, no.~6, pp. 1--35, 2021.

\bibitem{barocas-hardt-narayanan}
S.~Barocas, M.~Hardt, and A.~Narayanan, \emph{Fairness and Machine
  Learning}.\hskip 1em plus 0.5em minus 0.4em\relax fairmlbook.org, 2019,
  \url{http://www.fairmlbook.org}.

\bibitem{zehlike2021fairness}
M.~Zehlike, K.~Yang, and J.~Stoyanovich, ``Fairness in ranking: A survey,''
  \emph{arXiv preprint arXiv:2103.14000}, 2021.

\bibitem{limmer2018optimization}
S.~Limmer and M.~Dietrich, ``Optimization of dynamic prices for electric
  vehicle charging considering fairness,'' in \emph{2018 IEEE Symposium Series
  on Computational Intelligence (SSCI)}.\hskip 1em plus 0.5em minus 0.4em\relax
  IEEE, 2018, pp. 2304--2311.

\bibitem{2004TheUnfairPrice}
X.~Lan, K.~B. Monroe, and J.~L. Cox, ``The price is unfair! a conceptual
  framework of price fairness perceptions,'' \emph{Journal of Marketing},
  vol.~68, no.~4, pp. 1--15, 2004.

\bibitem{8285280}
S.~Limmer and T.~Rodemann, ``Multi-objective optimization of plug-in electric
  vehicle charging prices,'' in \emph{2017 IEEE Symposium Series on
  Computational Intelligence (SSCI)}, 2017, pp. 1--8.

\bibitem{FairCAS2007}
M.~Saini and S.~Rao, ``Fairness in combinatorial auctioning systems,''
  \emph{AAAI Spring Symposium - Technical Report}, pp. 61--67, 2007.

\bibitem{fairRA2021}
D.~Kumar, G.~Baranwal, and D.~P. Vidyarthi, ``Fair resource allocation policies
  in reverse auction-based cloud market,'' \emph{SN Computer Science}, vol.~2,
  p. 483, 2021.

\bibitem{ghodsi2011dominant}
A.~Ghodsi, M.~Zaharia, B.~Hindman, A.~Konwinski, S.~Shenker, and I.~Stoica,
  ``Dominant resource fairness: {F}air allocation of multiple resource types,''
  in \emph{8th USENIX Symposium on Networked Systems Design and Implementation
  (NSDI 11)}, 2011, pp. 323--336.

\bibitem{6919321}
W.~Wang, B.~Liang, and B.~Li, ``Multi-resource fair allocation in heterogeneous
  cloud computing systems,'' \emph{IEEE Transactions on Parallel and
  Distributed Systems}, vol.~26, no.~10, pp. 2822--2835, 2015.

\bibitem{6211077}
L.~Xu, Z.~Zeng, and X.~Ye, ``Multi-objective optimization based virtual
  resource allocation strategy for cloud computing,'' in \emph{2012 IEEE/ACIS
  11th International Conference on Computer and Information Science}.\hskip 1em
  plus 0.5em minus 0.4em\relax IEEE, 2012, pp. 56--61.

\bibitem{JIANG2016239}
\BIBentryALTinterwordspacing
H.~Jiang, J.~Yi, S.~Chen, and X.~Zhu, ``A multi-objective algorithm for task
  scheduling and resource allocation in cloud-based disassembly,''
  \emph{Journal of Manufacturing Systems}, vol.~41, pp. 239--255, 2016.
  [Online]. Available:
  \url{https://www.sciencedirect.com/science/article/pii/S0278612516300656}
\BIBentrySTDinterwordspacing

\bibitem{xu2021federatedEA}
J.~Xu, Y.~Jin, W.~Du, and S.~Gu, ``A federated data-driven evolutionary
  algorithm,'' \emph{Knowledge-Based Systems}, vol. 233, p. 107532, 2021.

\bibitem{xu2021federatedMOEA}
J.~Xu, Y.~Jin, and W.~Du, ``A federated data-driven evolutionary algorithm for
  expensive multi-/many-objective optimization,'' \emph{Complex \& Intelligent
  Systems}, vol.~7, no.~6, pp. 3093--3109, 2021.

\bibitem{hu2015face}
G.~Hu, Y.~Yang, D.~Yi, J.~Kittler, W.~Christmas, S.~Z. Li, and T.~Hospedales,
  ``When face recognition meets with deep learning: an evaluation of
  convolutional neural networks for face recognition,'' in \emph{Proceedings of
  the IEEE International Conference on Computer Vision Workshops}, 2015, pp.
  142--150.

\bibitem{davidson2010youtube}
J.~Davidson, B.~Liebald, J.~Liu, P.~Nandy, T.~Van~Vleet, U.~Gargi, S.~Gupta,
  Y.~He, M.~Lambert, B.~Livingston \emph{et~al.}, ``The youtube video
  recommendation system,'' in \emph{Proceedings of the Fourth ACM Conference on
  Recommender Systems}, 2010, pp. 293--296.

\bibitem{8456559}
Y.~Jin, H.~Wang, T.~Chugh, D.~Guo, and K.~Miettinen, ``Data-driven evolutionary
  optimization: An overview and case studies,'' \emph{IEEE Transactions on
  Evolutionary Computation}, vol.~23, no.~3, pp. 442--458, 2019.

\bibitem{jin2021data}
Y.~Jin, H.~Wang, and C.~Sun, \emph{Data-driven evolutionary
  optimization}.\hskip 1em plus 0.5em minus 0.4em\relax Springer International
  Publishing, 2021.

\bibitem{zheng2017decomposition}
J.~Zheng, G.~Yu, Q.~Zhu, X.~Li, and J.~Zou, ``On decomposition methods in
  interactive user-preference based optimization,'' \emph{Applied Soft
  Computing}, vol.~52, pp. 952--973, 2017.

\bibitem{8412189}
B.~Xin, L.~Chen, J.~Chen, H.~Ishibuchi, K.~Hirota, and B.~Liu, ``Interactive
  multiobjective optimization: A review of the state-of-the-art,'' \emph{IEEE
  Access}, vol.~6, pp. 41\,256--41\,279, 2018.

\bibitem{8955593}
W.~Zhang, X.~Tang, and J.~Wang, ``On fairness-aware learning for
  non-discriminative decision-making,'' in \emph{2019 International Conference
  on Data Mining Workshops (ICDMW)}.\hskip 1em plus 0.5em minus 0.4em\relax
  IEEE, 2019, pp. 1072--1079.

\bibitem{8710313}
M.~K. Tomczyk and M.~Kadziński, ``Decomposition-based interactive evolutionary
  algorithm for multiple objective optimization,'' \emph{IEEE Transactions on
  Evolutionary Computation}, vol.~24, no.~2, pp. 320--334, 2020.

\bibitem{2011Review}
P.~Agarwal, M.~Sahai, V.~Mishra, M.~Bag, and V.~Singh, ``A review of
  multi-criteria decision making techniques for supplier evaluation and
  selection,'' \emph{International Journal of Industrial Engineering
  Computations}, vol.~2, no.~4, pp. 801--810, 2011.

\bibitem{BECHIKH2015141}
S.~Bechikh, M.~Kessentini, L.~B. Said, and K.~Ghédira, ``Preference
  incorporation in evolutionary multiobjective optimization,'' \emph{Advances
  in Computers}, vol.~98, pp. 141--207, 2015.

\bibitem{J2008Multiobjective}
J.~Branke, K.~Deb, K.~Miettinen, and R.~Sowiński, \emph{Multiobjective
  Optimization: Interactive and Evolutionary Approaches}.\hskip 1em plus 0.5em
  minus 0.4em\relax Berlin, Germany: Springer, 2008.

\bibitem{YU2017689}
\BIBentryALTinterwordspacing
G.~Yu, R.~Shen, J.~Zheng, M.~Li, J.~Zou, and Y.~Liu, ``Binary search based
  boundary elimination selection in many-objective evolutionary optimization,''
  \emph{Applied Soft Computing}, vol.~60, pp. 689--705, 2017. [Online].
  Available:
  \url{https://www.sciencedirect.com/science/article/pii/S1568494617304416}
\BIBentrySTDinterwordspacing

\bibitem{9332241}
L.~Ma, M.~Huang, S.~Yang, R.~Wang, and X.~Wang, ``An adaptive localized
  decision variable analysis approach to large-scale multiobjective and
  many-objective optimization,'' \emph{IEEE Transactions on Cybernetics}, pp.
  1--13, 2021.

\bibitem{9484680}
L.~Ma, N.~Li, Y.~Guo, X.~Wang, S.~Yang, M.~Huang, and H.~Zhang, ``Learning to
  optimize: Reference vector reinforcement learning adaption to constrained
  many-objective optimization of industrial copper burdening system,''
  \emph{IEEE Transactions on Cybernetics}, pp. 1--14, 2021.

\bibitem{2003Performance}
V.~Khare, X.~Yao, and K.~Deb, ``Performance scaling of multi-objective
  evolutionary algorithms,'' in \emph{International Conference on Evolutionary
  Multi-Criterion Optimization}.\hskip 1em plus 0.5em minus 0.4em\relax
  Springer, Berlin, Heidelberg, 2003, p. 367–390.

\bibitem{ishibuchi2008evolutionary}
H.~Ishibuchi, N.~Tsukamoto, and Y.~Nojima, ``Evolutionary many-objective
  optimization: A short review,'' in \emph{Proceedings of the 2008 IEEE
  Congress on Evolutionary Computation}.\hskip 1em plus 0.5em minus 0.4em\relax
  Hong Kong, China: IEEE, 2008, pp. 2419--2426.

\bibitem{2002Axiomatization}
R.~Slowinski, S.~Greco, and B.~Matarazzo, ``Axiomatization of utility,
  outranking and decision rule preference models for multiple-criteria
  classification problems under partial inconsistency with the dominance
  principle,'' \emph{Control \& Cybernetics}, vol.~31, no.~4, pp. 1005--1035,
  2002.

\bibitem{4358754}
Q.~Zhang and H.~Li, ``{MOEA/D}: A multiobjective evolutionary algorithm based
  on decomposition,'' \emph{IEEE Transactions on Evolutionary Computation},
  vol.~11, no.~6, pp. 712--731, 2007.

\bibitem{7386636}
R.~Cheng, Y.~Jin, M.~Olhofer, and B.~Sendhoff, ``A reference vector guided
  evolutionary algorithm for many-objective optimization,'' \emph{IEEE
  Transactions on Evolutionary Computation}, vol.~20, no.~5, pp. 773--791,
  2016.

\bibitem{deb2006reference}
K.~Deb and J.~Sundar, ``Reference point based multi-objective optimization
  using evolutionary algorithms,'' in \emph{Proceedings of the 8th Annual
  Conference on Genetic and Evolutionary Computation}, ACM.\hskip 1em plus
  0.5em minus 0.4em\relax New York, NY, USA: Association for Computing
  Machinery, 2006, pp. 635--642.

\bibitem{knee9209056}
Z.~He, G.~G. Yen, and J.~Ding, ``Knee-based decision making and visualization
  in many-objective optimization,'' \emph{IEEE Transactions on Evolutionary
  Computation}, vol.~25, no.~2, pp. 292--306, 2021.

\bibitem{ruiz2008additive}
F.~Ruiz, M.~Luque, F.~Miguel, and M.~del Mar~Mu{\~n}oz, ``An additive
  achievement scalarizing function for multiobjective programming problems,''
  \emph{European Journal of Operational Research}, vol. 188, no.~3, pp.
  683--694, 2008.

\bibitem{Bhattacharjee2017Bridging}
K.~S. Bhattacharjee, H.~K. Singh, M.~Ryan, and T.~Ray, ``Bridging the gap:
  Many-objective optimization and informed decision-making,'' \emph{IEEE
  Transactions on Evolutionary Computation}, vol.~21, no.~5, pp. 813--820,
  2017.

\bibitem{7465803}
W.~{Chiu}, G.~G. {Yen}, and T.~{Juan}, ``Minimum manhattan distance approach to
  multiple criteria decision making in multiobjective optimization problems,''
  \emph{IEEE Transactions on Evolutionary Computation}, vol.~20, no.~6, pp.
  972--985, 2016.

\bibitem{8933061}
K.~{Zhang}, G.~G. {Yen}, and Z.~{He}, ``Evolutionary algorithm for knee-based
  multiple criteria decision making,'' \emph{IEEE Transactions on Cybernetics},
  vol.~51, no.~2, pp. 722--735, 2021.

\bibitem{2009gdominance}
J.~Molina, L.~V. Santana, A.~Hernández-Díaz, C.~Coello, and R.~Caballero,
  ``g-dominance: Reference point based dominance for multiobjective
  metaheuristics,'' \emph{European Journal of Operational Research}, vol. 197,
  no.~2, pp. 685--692, 2009.

\bibitem{cone2010}
A.~Sinha, P.~Korhonen, J.~Wallenius, and K.~Deb, ``An interactive evolutionary
  multi-objective optimization method based on polyhedral cones,'' in
  \emph{Learning and Intelligent Optimization}, C.~Blum and R.~Battiti,
  Eds.\hskip 1em plus 0.5em minus 0.4em\relax Berlin, Heidelberg: Springer
  Berlin Heidelberg, 2010, pp. 318--332.

\bibitem{2010rdominance}
L.~B. Said, S.~Bechikh, and K.~Ghedira, ``The r-dominance: A new dominance
  relation for interactive evolutionary multicriteria decision making,''
  \emph{IEEE Transactions on Evolutionary Computation}, vol.~14, no.~5, pp.
  801--818, 2010.

\bibitem{8552449}
J.~Yi, J.~Bai, H.~He, J.~Peng, and D.~Tang, ``ar-{MOEA}: A novel
  preference-based dominance relation for evolutionary multiobjective
  optimization,'' \emph{IEEE Transactions on Evolutionary Computation},
  vol.~23, no.~5, pp. 788--802, 2019.

\bibitem{G3321930}
G.~Yu, Y.~Jin, and M.~Olhofer, ``An a priori knee identification
  multi-objective evolutionary algorithm based on $\alpha$-dominance,'' in
  \emph{Proceedings of the Genetic and Evolutionary Computation Conference
  Companion}, ser. GECCO '19.\hskip 1em plus 0.5em minus 0.4em\relax New York,
  NY, USA: Association for Computing Machinery, 2019, pp. 241--242.

\bibitem{9139367}
G.~{Yu}, Y.~{Jin}, and M.~{Olhofer}, ``A multiobjective evolutionary algorithm
  for finding knee regions using two localized dominance relationships,''
  \emph{IEEE Transactions on Evolutionary Computation}, vol.~25, no.~1, pp.
  145--158, 2021.

\bibitem{lopez2009study}
A.~L{\'o}pez~Jaimes and C.~A. Coello~Coello, ``Study of preference relations in
  many-objective optimization,'' in \emph{Proceedings of the 11th Annual
  Conference on Genetic and Evolutionary Computation}.\hskip 1em plus 0.5em
  minus 0.4em\relax Montréal, Canada: ACM, 2009, pp. 611--618.

\bibitem{2001Rough}
S.~Greco, B.~Matarazzo, and R.~Slowinski, ``Rough sets theory for multicriteria
  decision analysis,'' \emph{European Journal of Operational Research}, vol.
  129, no.~1, pp. 1--47, 2001.

\bibitem{5585982}
S.~Greco, B.~Matarazzo, and R.~Slowiński, ``Interactive evolutionary
  multiobjective optimization using dominance-based rough set approach,'' in
  \emph{IEEE Congress on Evolutionary Computation}.\hskip 1em plus 0.5em minus
  0.4em\relax IEEE, 2010, pp. 1--8.

\bibitem{Greco2008}
\BIBentryALTinterwordspacing
S.~Greco, B.~Matarazzo, and R.~S{\l}owi{\'{n}}ski, \emph{Dominance-Based Rough
  Set Approach to Interactive Multiobjective Optimization}.\hskip 1em plus
  0.5em minus 0.4em\relax Berlin, Heidelberg: Springer Berlin Heidelberg, 2008,
  pp. 121--155. [Online]. Available:
  \url{https://doi.org/10.1007/978-3-540-88908-3_5}
\BIBentrySTDinterwordspacing

\bibitem{yaochu2002fuzzy}
Y.~Jin and B.~Sendhoff, ``Incorporation of fuzzy preferences into evolutionary
  multiobjective optimization,'' in \emph{Proceedings of the 4th Annual
  Conference on Genetic and Evolutionary Computation}, ser. GECCO'02.\hskip 1em
  plus 0.5em minus 0.4em\relax San Francisco, CA, USA: Morgan Kaufmann
  Publishers Inc., 2002, p. 683.

\bibitem{2012Expressing}
A.~Hadjali, A.~Mokhtari, and O.~Pivert, ``Expressing and processing complex
  preferences in route planning queries: Towards a fuzzy-set-based approach,''
  \emph{Fuzzy Sets \& Systems}, vol. 196, pp. 82--104, 2012.

\bibitem{2011Relational}
R.~I. Brafman, ``Relational preference rules for control,'' \emph{Artificial
  Intelligence}, vol. 175, no. 7–8, pp. 1180--1193, 2011.

\bibitem{2010Brain}
R.~Battiti and A.~Passerini, ``Brain–computer evolutionary multiobjective
  optimization: A genetic algorithm adapting to the decision maker,''
  \emph{IEEE Transactions on Evolutionary Computation}, vol.~14, no.~5, pp.
  671--687, 2010.

\bibitem{yu2016decomposing}
G.~Yu, J.~Zheng, R.~Shen, and M.~Li, ``Decomposing the user-preference in
  multiobjective optimization,'' \emph{Soft Computing}, vol.~20, no.~10, pp.
  4005--4021, 2016.

\bibitem{Jianjie2017A}
J.~Hu, G.~Yu, J.~Zheng, and J.~Zou, ``A preference-based multi-objective
  evolutionary algorithm using preference selection radius,'' \emph{Soft
  Computing}, vol.~21, pp. 5025--5051, 2017.

\bibitem{5585741}
M.~K{\"o}ksalan and {\.I}.~Karahan, ``An interactive territory defining
  evolutionary algorithm: i{TDEA},'' \emph{IEEE Transactions on Evolutionary
  Computation}, vol.~14, no.~5, pp. 702--722, 2010.

\bibitem{7572016}
R.~Wang, Z.~Zhou, H.~Ishibuchi, T.~Liao, and T.~Zhang, ``Localized weighted sum
  method for many-objective optimization,'' \emph{IEEE Transactions on
  Evolutionary Computation}, vol.~22, no.~1, pp. 3--18, 2018.

\bibitem{Branke2005}
\BIBentryALTinterwordspacing
J.~Branke and K.~Deb, \emph{Integrating User Preferences into Evolutionary
  Multi-Objective Optimization}.\hskip 1em plus 0.5em minus 0.4em\relax Berlin,
  Heidelberg: Springer Berlin Heidelberg, 2005, pp. 461--477. [Online].
  Available: \url{https://doi.org/10.1007/978-3-540-44511-1_21}
\BIBentrySTDinterwordspacing

\bibitem{zitzler2008spam}
E.~Zitzler, L.~Thiele, and J.~Bader, ``{SPAM}: Set preference algorithm for
  multiobjective optimization,'' in \emph{Proceedings of the 10th International
  Conference on Parallel Problem Solving from Nature}.\hskip 1em plus 0.5em
  minus 0.4em\relax Dortmund, Germany: Springer, Berlin, Heidelberg, 2008, pp.
  847--858.

\bibitem{saaty2008decision}
T.~L. Saaty, ``Decision making with the analytic hierarchy process,''
  \emph{International Journal of Services Sciences}, vol.~1, no.~1, pp. 83--98,
  2008.

\bibitem{6600851}
K.~Deb and H.~Jain, ``An evolutionary many-objective optimization algorithm
  using reference-point-based nondominated sorting approach, {P}art {I}:
  Solving problems with box constraints,'' \emph{IEEE Transactions on
  Evolutionary Computation}, vol.~18, no.~4, pp. 577--601, 2014.

\bibitem{Miettinen2008}
\BIBentryALTinterwordspacing
K.~Miettinen, F.~Ruiz, and A.~P. Wierzbicki, \emph{Introduction to
  Multiobjective Optimization: Interactive Approaches}.\hskip 1em plus 0.5em
  minus 0.4em\relax Berlin, Heidelberg: Springer Berlin Heidelberg, 2008, pp.
  27--57. [Online]. Available:
  \url{https://doi.org/10.1007/978-3-540-88908-3_2}
\BIBentrySTDinterwordspacing

\bibitem{siskos1986use}
J.~Siskos, J.~Lombard, and A.~Oudiz, ``The use of multicriteria outranking
  methods in the comparison of control options against a chemical pollutant,''
  \emph{Journal of the Operational Research Society}, vol.~37, no.~4, pp.
  357--371, 1986.

\bibitem{waegeman2011era}
W.~Waegeman and B.~De~Baets, ``On the {ERA} ranking representability of
  pairwise bipartite ranking functions,'' \emph{Artificial Intelligence}, vol.
  175, no. 7-8, pp. 1223--1250, 2011.

\bibitem{2008Ordinal}
S.~Greco, V.~Mousseau, and R.~Slowinski, ``Ordinal regression revisited:
  Multiple criteria ranking using a set of additive value functions,''
  \emph{European Journal of Operational Research}, vol. 191, no.~2, pp.
  416--436, 2008.

\bibitem{8638825}
G.~{Yu}, Y.~{Jin}, and M.~{Olhofer}, ``Benchmark problems and performance
  indicators for search of knee points in multiobjective optimization,''
  \emph{IEEE Transactions on Cybernetics}, vol.~50, no.~8, pp. 3531--3544,
  2020.

\bibitem{2011Understanding}
K.~Deb and S.~Gupta, ``Understanding knee points in bicriteria problems and
  their implications as preferred solution principles,'' \emph{Engineering
  Optimization}, vol.~43, no.~11, pp. 1175--1204, 2011.

\bibitem{1999On}
I.~Das, ``On characterizing the ``knee'' of the {P}areto curve based on
  normal-boundary intersection,'' \emph{Structural Optimization}, vol.~18, no.
  2-3, pp. 107--115, 1999.

\bibitem{caton2020fairness}
S.~Caton and C.~Haas, ``Fairness in machine learning: A survey,'' \emph{arXiv
  preprint arXiv:2010.04053}, 2020.

\bibitem{quy2021survey}
T.~L. Quy, A.~Roy, V.~Iosifidis, and E.~Ntoutsi, ``A survey on datasets for
  fairness-aware machine learning,'' \emph{arXiv preprint arXiv:2110.00530},
  2021.

\bibitem{olteanu2019social}
\BIBentryALTinterwordspacing
A.~Olteanu, C.~Castillo, F.~Diaz, and E.~K{\i}c{\i}man, ``Social data: Biases,
  methodological pitfalls, and ethical boundaries,'' \emph{Frontiers in Big
  Data}, vol.~2, 2019. [Online]. Available:
  \url{https://doi.org/10.3389/fdata.2019.00013}
\BIBentrySTDinterwordspacing

\bibitem{zhang2018mitigating}
B.~H. Zhang, B.~Lemoine, and M.~Mitchell, ``Mitigating unwanted biases with
  adversarial learning,'' in \emph{Proceedings of the 2018 AAAI/ACM Conference
  on AI, Ethics, and Society}, 2018, pp. 335--340.

\bibitem{suresh2019framework}
H.~Suresh and J.~V. Guttag, ``A framework for understanding unintended
  consequences of machine learning,'' \emph{arXiv preprint arXiv:1901.10002},
  vol.~2, 2019.

\bibitem{clarke2005phantom}
K.~A. Clarke, ``The phantom menace: Omitted variable bias in econometric
  research,'' \emph{Conflict Management and Peace Science}, vol.~22, no.~4, pp.
  341--352, 2005.

\bibitem{buolamwini2018gender}
J.~Buolamwini and T.~Gebru, ``Gender shades: Intersectional accuracy
  disparities in commercial gender classification,'' in \emph{Conference on
  Fairness, Accountability and Transparency}.\hskip 1em plus 0.5em minus
  0.4em\relax PMLR, 2018, pp. 77--91.

\bibitem{danks2017algorithmic}
D.~Danks and A.~J. London, ``Algorithmic bias in autonomous systems.'' in
  \emph{IJCAI}, vol.~17, 2017, pp. 4691--4697.

\bibitem{ciampaglia2018algorithmic}
G.~L. Ciampaglia, A.~Nematzadeh, F.~Menczer, and A.~Flammini, ``How algorithmic
  popularity bias hinders or promotes quality,'' \emph{Scientific Reports},
  vol.~8, no.~1, pp. 1--7, 2018.

\bibitem{nguyen2013old}
D.~Nguyen, R.~Gravel, D.~Trieschnigg, and T.~Meder, ``"how old do you think {I}
  am?" a study of language and age in {T}witter,'' in \emph{Proceedings of the
  International AAAI Conference on Web and Social Media}, vol.~7, no.~1.\hskip
  1em plus 0.5em minus 0.4em\relax AAAI Press, Palo Alto, California USA, 2013,
  pp. 439--448.

\bibitem{dwork2012fairness}
C.~Dwork, M.~Hardt, T.~Pitassi, O.~Reingold, and R.~Zemel, ``Fairness through
  awareness,'' in \emph{Proceedings of the 3rd Innovations in Theoretical
  Computer Science Conference}.\hskip 1em plus 0.5em minus 0.4em\relax New
  York, NY, USA: Association for Computing Machinery, 2012, pp. 214--226.

\bibitem{kamishima2012fairness}
T.~Kamishima, S.~Akaho, H.~Asoh, and J.~Sakuma, ``Fairness-aware classifier
  with prejudice remover regularizer,'' in \emph{Joint European Conference on
  Machine Learning and Knowledge Discovery in Databases}, Springer.\hskip 1em
  plus 0.5em minus 0.4em\relax Berlin, Heidelberg: Springer-Verlag, 2012, pp.
  35--50.

\bibitem{zemel2013learning}
R.~Zemel, Y.~Wu, K.~Swersky, T.~Pitassi, and C.~Dwork, ``Learning fair
  representations,'' in \emph{International Conference on Machine Learning},
  PMLR.\hskip 1em plus 0.5em minus 0.4em\relax JMLR.org, 2013, pp. 325--333.

\bibitem{hardt2016equality}
M.~Hardt, E.~Price, and N.~Srebro, ``Equality of opportunity in supervised
  learning,'' in \emph{Proceedings of the 30th International Conference on
  Neural Information Processing Systems}.\hskip 1em plus 0.5em minus
  0.4em\relax Red Hook, NY, USA: Curran Associates Inc., 2016, pp. 3323--3331.

\bibitem{verma2018fairness}
S.~Verma and J.~Rubin, ``Fairness definitions explained,'' in \emph{2018
  IEEE/ACM International Workshop on Software Fairness (Fairware)}.\hskip 1em
  plus 0.5em minus 0.4em\relax IEEE, 2018, pp. 1--7.

\bibitem{corbett2017algorithmic}
S.~Corbett-Davies, E.~Pierson, A.~Feller, S.~Goel, and A.~Huq, ``Algorithmic
  decision making and the cost of fairness,'' in \emph{Proceedings of the 23rd
  ACM SIGKDD International Conference on Knowledge Discovery and Data
  Mining}.\hskip 1em plus 0.5em minus 0.4em\relax New York, NY, USA:
  Association for Computing Machinery, 2017, pp. 797--806.

\bibitem{kearns2018preventing}
M.~Kearns, S.~Neel, A.~Roth, and Z.~S. Wu, ``Preventing fairness
  gerrymandering: Auditing and learning for subgroup fairness,'' in
  \emph{Proceedings of the 35th International Conference on Machine Learning},
  vol.~80.\hskip 1em plus 0.5em minus 0.4em\relax PMLR, 2018, pp. 2564--2572.

\bibitem{kearns2019empirical}
------, ``An empirical study of rich subgroup fairness for machine learning,''
  in \emph{Proceedings of the Conference on Fairness, Accountability, and
  Transparency}.\hskip 1em plus 0.5em minus 0.4em\relax New York, NY, USA:
  Association for Computing Machinery, 2019, pp. 100--109.

\bibitem{grgic2016case}
N.~Grgic-Hlaca, M.~B. Zafar, K.~P. Gummadi, and A.~Weller, ``The case for
  process fairness in learning: Feature selection for fair decision making,''
  in \emph{NIPS Symposium on Machine Learning and the Law}, vol.~8, 2016.

\bibitem{kusner2017counterfactual}
M.~J. Kusner, J.~Loftus, C.~Russell, and R.~Silva, ``Counterfactual fairness,''
  \emph{Advances in Neural Information Processing Systems}, vol.~30, no.
  4066--4076, 2017.

\bibitem{calmon2017optimized}
F.~Calmon, D.~Wei, B.~Vinzamuri, K.~Natesan~Ramamurthy, and K.~R. Varshney,
  ``Optimized pre-processing for discrimination prevention,'' \emph{Advances in
  Neural Information Processing Systems}, vol.~30, pp. 3992--4001, 2017.

\bibitem{creager2019flexibly}
E.~Creager, D.~Madras, J.-H. Jacobsen, M.~Weis, K.~Swersky, T.~Pitassi, and
  R.~Zemel, ``Flexibly fair representation learning by disentanglement,'' in
  \emph{International Conference on Machine Learning}.\hskip 1em plus 0.5em
  minus 0.4em\relax PMLR, 2019, pp. 1436--1445.

\bibitem{kamiran2012data}
F.~Kamiran and T.~Calders, ``Data preprocessing techniques for classification
  without discrimination,'' \emph{Knowledge and Information Systems}, vol.~33,
  no.~1, pp. 1--33, 2012.

\bibitem{feldman2015certifying}
M.~Feldman, S.~A. Friedler, J.~Moeller, C.~Scheidegger, and
  S.~Venkatasubramanian, ``Certifying and removing disparate impact,'' in
  \emph{proceedings of the 21th ACM SIGKDD International Conference on
  Knowledge Discovery and Data Mining}.\hskip 1em plus 0.5em minus 0.4em\relax
  New York, NY, USA: Association for Computing Machinery, 2015, pp. 259--268.

\bibitem{bellamy2019ai}
R.~K.~E. Bellamy, K.~Dey, M.~Hind, S.~C. Hoffman, S.~Houde, K.~Kannan,
  P.~Lohia, J.~Martino, S.~Mehta, A.~Mojsilović, S.~Nagar, K.~N. Ramamurthy,
  J.~Richards, D.~Saha, P.~Sattigeri, M.~Singh, K.~R. Varshney, and Y.~Zhang,
  ``Ai fairness 360: An extensible toolkit for detecting and mitigating
  algorithmic bias,'' \emph{IBM Journal of Research and Development}, vol.~63,
  no. 4/5, pp. 4:1--4:15, 2019.

\bibitem{chawla2002smote}
N.~V. Chawla, K.~W. Bowyer, L.~O. Hall, and W.~P. Kegelmeyer, ``Smote:
  synthetic minority over-sampling technique,'' \emph{Journal of Artificial
  Intelligence Research}, vol.~16, pp. 321--357, 2002.

\bibitem{dixon2018measuring}
L.~Dixon, J.~Li, J.~Sorensen, N.~Thain, and L.~Vasserman, ``Measuring and
  mitigating unintended bias in text classification,'' in \emph{Proceedings of
  the 2018 AAAI/ACM Conference on AI, Ethics, and Society}.\hskip 1em plus
  0.5em minus 0.4em\relax New York, NY, USA: Association for Computing
  Machinery, 2018, pp. 67--73.

\bibitem{d2017conscientious}
B.~d'Alessandro, C.~O'Neil, and T.~LaGatta, ``Conscientious classification: A
  data scientist's guide to discrimination-aware classification,'' \emph{Big
  Data}, vol.~5, no.~2, pp. 120--134, 2017.

\bibitem{cotter2019two}
A.~Cotter, H.~Jiang, and K.~Sridharan, ``Two-player games for efficient
  non-convex constrained optimization,'' in \emph{Algorithmic Learning
  Theory}.\hskip 1em plus 0.5em minus 0.4em\relax PMLR, 2019, pp. 300--332.

\bibitem{liu2022accuracy}
S.~Liu and L.~N. Vicente, ``Accuracy and fairness trade-offs in machine
  learning: A stochastic multi-objective approach,'' \emph{Computational
  Management Science}, pp. 1--25, 2022.

\bibitem{fish2016confidence}
B.~Fish, J.~Kun, and {\'A}.~D. Lelkes, ``A confidence-based approach for
  balancing fairness and accuracy,'' in \emph{Proceedings of the 2016 SIAM
  International Conference on Data Mining}.\hskip 1em plus 0.5em minus
  0.4em\relax SIAM, 2016, pp. 144--152.

\bibitem{pleiss2017fairness}
G.~Pleiss, M.~Raghavan, F.~Wu, J.~Kleinberg, and K.~Q. Weinberger, ``On
  fairness and calibration,'' \emph{Advances in Neural Information Processing
  Systems}, p. 5684–5693, 2017.

\bibitem{kamiran2012decision}
F.~Kamiran, A.~Karim, and X.~Zhang, ``Decision theory for discrimination-aware
  classification,'' in \emph{2012 IEEE 12th International Conference on Data
  Mining}.\hskip 1em plus 0.5em minus 0.4em\relax IEEE, 2012, pp. 924--929.

\bibitem{1993TradeOffs}
L.~D. Ordoñez and B.~A. Mellers, ``Trade-offs in fairness and preference
  judgments,'' \emph{Psychological Perspectives on Justice: Theory and
  Applications}, pp. 138--154, 1993.

\bibitem{avin2015homophily}
C.~Avin, B.~Keller, Z.~Lotker, C.~Mathieu, D.~Peleg, and Y.-A. Pignolet,
  ``Homophily and the glass ceiling effect in social networks,'' in
  \emph{Proceedings of the 2015 Conference on Innovations in Theoretical
  Computer Science}.\hskip 1em plus 0.5em minus 0.4em\relax New York, NY, USA:
  Association for Computing Machinery, 2015, pp. 41--50.

\bibitem{kim2019preference}
M.~P. Kim, A.~Korolova, G.~N. Rothblum, and G.~Yona, ``Preference-informed
  fairness,'' \emph{arXiv preprint arXiv:1904.01793}, 2019.

\bibitem{finocchiaro2021bridging}
J.~Finocchiaro, R.~Maio, F.~Monachou, G.~K. Patro, M.~Raghavan, A.-A. Stoica,
  and S.~Tsirtsis, ``Bridging machine learning and mechanism design towards
  algorithmic fairness,'' in \emph{Proceedings of the 2021 ACM Conference on
  Fairness, Accountability, and Transparency}.\hskip 1em plus 0.5em minus
  0.4em\relax New York, NY, USA: Association for Computing Machinery, 2021, pp.
  489--503.

\bibitem{do2021online}
V.~Do, S.~Corbett-Davies, J.~Atif, and N.~Usunier, ``Online certification of
  preference-based fairness for personalized recommender systems,'' \emph{arXiv
  preprint arXiv:2104.14527}, 2021.

\bibitem{zafar2017parity}
M.~B. Zafar, I.~Valera, M.~Rodriguez, K.~Gummadi, and A.~Weller, ``From parity
  to preference-based notions of fairness in classification,'' \emph{Advances
  in Neural Information Processing Systems}, vol.~30, pp. 229--239, 2017.

\bibitem{liu2019personalized}
W.~Liu, J.~Guo, N.~Sonboli, R.~Burke, and S.~Zhang, ``Personalized
  fairness-aware re-ranking for microlending,'' in \emph{Proceedings of the
  13th ACM Conference on Recommender Systems}.\hskip 1em plus 0.5em minus
  0.4em\relax New York, NY, USA: Association for Computing Machinery, 2019, pp.
  467--471.

\bibitem{saxena2019fairness}
N.~A. Saxena, K.~Huang, E.~DeFilippis, G.~Radanovic, D.~C. Parkes, and Y.~Liu,
  ``How do fairness definitions fare? examining public attitudes towards
  algorithmic definitions of fairness,'' in \emph{Proceedings of the 2019
  AAAI/ACM Conference on AI, Ethics, and Society}.\hskip 1em plus 0.5em minus
  0.4em\relax New York, NY, USA: Association for Computing Machinery, 2019, pp.
  99--106.

\bibitem{menon2017cost}
A.~K. Menon and R.~C. Williamson, ``The cost of fairness in classification,''
  \emph{arXiv preprint arXiv:1705.09055}, 2017.

\bibitem{valera2018enhancing}
I.~Valera, A.~Singla, and M.~Gomez~Rodriguez, ``Enhancing the accuracy and
  fairness of human decision making,'' \emph{Advances in Neural Information
  Processing Systems}, vol.~31, pp. 1769--1778, 2018.

\bibitem{laumanns2004running}
L.~Marco, T.~Lothar, and Z.~Eckart, ``Running time analysis of multiobjective
  evolutionary algorithms on pseudo-boolean functions,'' \emph{IEEE
  Transactions on Evolutionary Computation}, vol.~8, no.~2, pp. 170--182, 2004.

\bibitem{friedrich2008runtime}
T.~Friedrich, C.~Horoba, and F.~Neumann, ``Runtime analyses for using fairness
  in evolutionary multi-objective optimization,'' in \emph{International
  Conference on Parallel Problem Solving from Nature}, Springer.\hskip 1em plus
  0.5em minus 0.4em\relax Springer, Berlin, Heidelberg, 2008, pp. 671--680.

\bibitem{friedrich2011illustration}
------, ``Illustration of fairness in evolutionary multi-objective
  optimization,'' \emph{Theoretical Computer Science}, vol. 412, no.~17, pp.
  1546--1556, 2011.

\bibitem{qian2016selection}
C.~Qian, K.~Tang, and Z.-H. Zhou, ``Selection hyper-heuristics can provably be
  helpful in evolutionary multi-objective optimization,'' in
  \emph{International Conference on Parallel Problem Solving from Nature},
  Springer.\hskip 1em plus 0.5em minus 0.4em\relax Springer, Cham, 2016, pp.
  835--846.

\bibitem{limmer2017multi}
S.~Limmer and T.~Rodemann, ``Multi-objective optimization of plug-in electric
  vehicle charging prices,'' in \emph{2017 IEEE Symposium Series on
  Computational Intelligence (SSCI)}.\hskip 1em plus 0.5em minus 0.4em\relax
  IEEE, 2017, pp. 1--8.

\bibitem{lunz2012influence}
B.~Lunz, Z.~Yan, J.~B. Gerschler, and D.~U. Sauer, ``Influence of plug-in
  hybrid electric vehicle charging strategies on charging and battery
  degradation costs,'' \emph{Energy Policy}, vol.~46, pp. 511--519, 2012.

\bibitem{bagnall2001next}
A.~J. Bagnall, V.~J. Rayward-Smith, and I.~M. Whittley, ``The next release
  problem,'' \emph{Information and Software Technology}, vol.~43, no.~14, pp.
  883--890, 2001.

\bibitem{dong2022multi}
S.~Dong, Y.~Xue, S.~Brinkkemper, and Y.-F. Li, ``Multi-objective integer
  programming approaches to next release problem—enhancing exact methods for
  finding whole pareto front,'' \emph{Information and Software Technology},
  vol. 147, p. 106825, 2022.

\bibitem{finkelstein2009search}
A.~Finkelstein, M.~Harman, S.~A. Mansouri, J.~Ren, and Y.~Zhang, ``A search
  based approach to fairness analysis in requirement assignments to aid
  negotiation, mediation and decision making,'' \emph{Requirements
  Engineering}, vol.~14, no.~4, pp. 231--245, 2009.

\bibitem{ghasemi2021multi}
M.~Ghasemi, K.~Bagherifard, H.~Parvin, S.~Nejatian, and K.-H. Pho,
  ``Multi-objective whale optimization algorithm and multi-objective grey wolf
  optimizer for solving next release problem with developing fairness and
  uncertainty quality indicators,'' \emph{Applied Intelligence}, vol.~51,
  no.~8, pp. 5358--5387, 2021.

\bibitem{rossi2005aggregating}
F.~Rossi, K.~B. Venable, and T.~Walsh, ``Aggregating preferences cannot be
  fair.'' \emph{Intelligenza Artificiale}, vol.~2, no.~1, pp. 30--38, 2005.

\bibitem{zhang2005fairness}
Z.~Zhang, H.~Xing, Z.~Wang, and Q.~Ni, ``The fairness of ranking procedure in
  pair-wise preference learning,'' in \emph{2005 International Conference on
  Natural Language Processing and Knowledge Engineering}.\hskip 1em plus 0.5em
  minus 0.4em\relax IEEE, 2005, pp. 780--784.

\bibitem{luss1999equitable}
H.~Luss, ``On equitable resource allocation problems: A lexicographic minimax
  approach,'' \emph{Operations Research}, vol.~47, no.~3, pp. 361--378, 1999.

\bibitem{bowerman1995multi}
R.~Bowerman, B.~Hall, and P.~Calamai, ``A multi-objective optimization approach
  to urban school bus routing: Formulation and solution method,''
  \emph{Transportation Research Part A: Policy and Practice}, vol.~29, no.~2,
  pp. 107--123, 1995.

\bibitem{purushothaman2021evolutionary}
K.~Purushothaman and V.~Nagarajan, ``Evolutionary multi-objective optimization
  algorithm for resource allocation using deep neural network in 5g multi-user
  massive mimo,'' \emph{International Journal of Electronics}, vol. 108, no.~7,
  pp. 1214--1233, 2021.

\bibitem{marsh1994equity}
M.~T. Marsh and D.~A. Schilling, ``Equity measurement in facility location
  analysis: A review and framework,'' \emph{European Journal of Operational
  Research}, vol.~74, no.~1, pp. 1--17, 1994.

\bibitem{jagtenberg2021introducing}
C.~J. Jagtenberg, M.~A. Vollebergh, O.~Uleberg, and J.~R{\o}islien,
  ``Introducing fairness in norwegian air ambulance base location planning,''
  \emph{Scandinavian Journal of Trauma, Resuscitation and Emergency Medicine},
  vol.~29, no.~1, pp. 1--10, 2021.

\bibitem{2002Linear}
M.~M. Kostreva and W.~Ogryczak, ``Linear optimization with multiple equitable
  criteria,'' \emph{RAIRO - Operations Research}, vol.~33, no.~3, pp. 275--297,
  1999.

\bibitem{kostreva2004equitable}
M.~M. Kostreva, W.~Ogryczak, and A.~Wierzbicki, ``Equitable aggregations and
  multiple criteria analysis,'' \emph{European Journal of Operational
  Research}, vol. 158, no.~2, pp. 362--377, 2004.

\bibitem{ogryczak2014fair}
W.~Ogryczak, ``Fair optimization--methodological foundations of fairness in
  network resource allocation,'' in \emph{2014 IEEE 38th International Computer
  Software and Applications Conference Workshops}.\hskip 1em plus 0.5em minus
  0.4em\relax IEEE, 2014, pp. 43--48.

\bibitem{koppen2010comparison}
M.~K{\"o}ppen, R.~Verschae, K.~Yoshida, and M.~Tsuru, ``Comparison of
  evolutionary multi-objective optimization algorithms for the utilization of
  fairness in network control,'' in \emph{2010 IEEE International Conference on
  Systems, Man and Cybernetics}.\hskip 1em plus 0.5em minus 0.4em\relax IEEE,
  2010, pp. 2647--2655.

\bibitem{le2005rate}
J.-Y. Le~Boudec, ``Rate adaptation, congestion control and fairness: A
  tutorial,'' \emph{Web page, November}, 2005.

\bibitem{blanco2022fairness}
V.~Blanco and R.~G{\'a}zquez, ``Fairness in maximal covering facility location
  problems,'' \emph{arXiv preprint arXiv:2204.06446}, 2022.

\bibitem{jagtenberg2020fairness}
C.~Jagtenberg and A.~Mason, ``Fairness in the ambulance location problem:
  maximizing the bernoulli-nash social welfare,'' \emph{Available at SSRN
  3536707}, 2020.

\bibitem{burlacu2020solving}
R.~Burlacu, B.~Gei{\ss}ler, and L.~Schewe, ``Solving mixed-integer nonlinear
  programmes using adaptively refined mixed-integer linear programmes,''
  \emph{Optimization Methods and Software}, vol.~35, no.~1, pp. 37--64, 2020.

\bibitem{habibi2016multi}
M.~Habibi~Davijani, M.~Banihabib, A.~Nadjafzadeh~Anvar, and S.~Hashemi,
  ``Multi-objective optimization model for the allocation of water resources in
  arid regions based on the maximization of socioeconomic efficiency,''
  \emph{Water Resources Management}, vol.~30, no.~3, pp. 927--946, 2016.

\bibitem{solgi2020multi}
M.~Solgi, O.~Bozorg-Haddad, and H.~A. Lo{\'a}iciga, ``A multi-objective
  optimization model for operation of intermittent water distribution
  networks,'' \emph{Water Supply}, vol.~20, no.~7, pp. 2630--2647, 2020.

\bibitem{fu2018water}
J.~Fu, P.-A. Zhong, F.~Zhu, J.~Chen, Y.-n. Wu, and B.~Xu, ``Water resources
  allocation in transboundary river based on asymmetric nash--harsanyi
  leader--follower game model,'' \emph{Water}, vol.~10, no.~3, p. 270, 2018.

\bibitem{fu2021comparison}
J.~Fu, P.-A. Zhong, B.~Xu, F.~Zhu, J.~Chen, and J.~Li, ``Comparison of
  transboundary water resources allocation models based on game theory and
  multi-objective optimization,'' \emph{Water}, vol.~13, no.~10, p. 1421, 2021.

\bibitem{tang2021new}
X.~Tang, Y.~He, P.~Qi, Z.~Chang, M.~Jiang, and Z.~Dai, ``A new multi-objective
  optimization model of water resources considering fairness and water shortage
  risk,'' \emph{Water}, vol.~13, no.~19, p. 2648, 2021.

\bibitem{cullis2007applying}
J.~Cullis and B.~Van~Koppen, \emph{Applying the Gini coefficient to measure
  inequality of water use in the Olifants river water management area, South
  Africa}.\hskip 1em plus 0.5em minus 0.4em\relax IWMI, 2007, vol. 113.

\bibitem{li2022evolutionary}
J.-Y. Li, Z.-H. Zhan, and J.~Zhang, ``Evolutionary computation for expensive
  optimization: A survey,'' \emph{Machine Intelligence Research}, vol.~19,
  no.~1, pp. 3--23, 2022.

\bibitem{jiang2021interpretable}
C.~Jiang, R.~Vinuesa, R.~Chen, J.~Mi, S.~Laima, and H.~Li, ``An interpretable
  framework of data-driven turbulence modeling using deep neural networks,''
  \emph{Physics of Fluids}, vol.~33, no.~5, p. 055133, 2021.

\bibitem{sun2022learning}
H.~Sun, W.~Pu, X.~Fu, T.-H. Chang, and M.~Hong, ``Learning to continuously
  optimize wireless resource in a dynamic environment: A bilevel optimization
  perspective,'' \emph{IEEE Transactions on Signal Processing}, 2022.

\bibitem{hong2014signal}
M.~Hong and Z.-Q. Luo, ``Signal processing and optimal resource allocation for
  the interference channel,'' in \emph{Academic Press Library in Signal
  Processing}.\hskip 1em plus 0.5em minus 0.4em\relax Elsevier, 2014, vol.~2,
  pp. 409--469.

\bibitem{kaur2022vaxequity}
N.~Kaur, J.~Hughes, and J.~Chen, ``Vaxequity: A data-driven risk assessment and
  optimization framework for equitable vaccine distribution,'' \emph{arXiv
  preprint arXiv:2201.07321}, 2022.

\bibitem{5461911}
T.~Lan, D.~Kao, M.~Chiang, and A.~Sabharwal, ``An axiomatic theory of fairness
  in network resource allocation,'' in \emph{2010 Proceedings IEEE INFOCOM},
  2010, pp. 1--9.

\bibitem{wang2022survey}
X.~Wang, Y.~Jin, S.~Schmitt, and M.~Olhofer, ``Recent advances in {B}ayesian
  optimization,'' \emph{arXiv:2206.03301}, 2022.

\bibitem{perrone2021fair}
V.~Perrone, M.~Donini, M.~B. Zafar, R.~Schmucker, K.~Kenthapadi, and
  C.~Archambeau, ``Fair bayesian optimization,'' in \emph{Proceedings of the
  2021 AAAI/ACM Conference on AI, Ethics, and Society}, 2021, pp. 854--863.

\bibitem{williams2006gaussian}
C.~K. Williams and C.~E. Rasmussen, \emph{Gaussian processes for machine
  learning}.\hskip 1em plus 0.5em minus 0.4em\relax MIT press Cambridge, MA,
  2005.

\bibitem{gardner2014bayesian}
J.~R. Gardner, M.~J. Kusner, Z.~E. Xu, K.~Q. Weinberger, and J.~P. Cunningham,
  ``Bayesian optimization with inequality constraints.'' in \emph{Proceedings
  of the 31th International Conference on Machine Learning}, vol. 2014, 2014,
  pp. 937--945.

\bibitem{gelbart2014bayesian}
M.~A. Gelbart, J.~Snoek, and R.~P. Adams, ``Bayesian optimization with unknown
  constraints,'' \emph{arXiv preprint arXiv:1403.5607}, 2014.

\bibitem{donini2018empirical}
M.~Donini, L.~Oneto, S.~Ben-David, J.~S. Shawe-Taylor, and M.~Pontil,
  ``Empirical risk minimization under fairness constraints,'' \emph{Advances in
  Neural Information Processing Systems}, vol.~31, 2018.

\bibitem{zhu2021real}
H.~Zhu and Y.~Jin, ``Real-time federated evolutionary neural architecture
  search,'' \emph{IEEE Transactions on Evolutionary Computation}, vol.~26,
  no.~2, pp. 364--378, 2022.

\bibitem{2021FromFL}
H.~Zhu, H.~Zhang, and Y.~Jin, ``From federated learning to federated neural
  architecture search: a survey,'' \emph{Complex \& Intelligent Systems},
  vol.~7, no.~2, pp. 639--657, 2021.

\bibitem{yue2021gifair}
X.~Yue, M.~Nouiehed, and R.~A. Kontar, ``{GIFAIR-FL}: An approach for group and
  individual fairness in federated learning,'' \emph{arXiv preprint
  arXiv:2108.02741}, 2021.

\bibitem{li2019fair}
T.~Li, M.~Sanjabi, A.~Beirami, and V.~Smith, ``Fair resource allocation in
  federated learning,'' \emph{arXiv preprint arXiv:1905.10497}, 2019.

\bibitem{lyu2020collaborative}
L.~Lyu, X.~Xu, Q.~Wang, and H.~Yu, ``Collaborative fairness in federated
  learning,'' in \emph{Federated Learning}.\hskip 1em plus 0.5em minus
  0.4em\relax Springer, 2020, pp. 189--204.

\bibitem{ntoutsi2020bias}
E.~Ntoutsi, P.~Fafalios, U.~Gadiraju, V.~Iosifidis, W.~Nejdl, M.-E. Vidal,
  S.~Ruggieri, F.~Turini, S.~Papadopoulos, E.~Krasanakis \emph{et~al.}, ``Bias
  in data-driven artificial intelligence systems—an introductory survey,''
  \emph{Wiley Interdisciplinary Reviews: Data Mining and Knowledge Discovery},
  vol.~10, no.~3, p. e1356, 2020.

\bibitem{2006Multi}
Y.~Jin, \emph{Multi-Objective Machine Learning}.\hskip 1em plus 0.5em minus
  0.4em\relax Multi-Objective Machine Learning, 2006.

\bibitem{gu2015ensemble}
S.~Gu, R.~Cheng, and Y.~Jin, ``Multi-objective ensemble generation,''
  \emph{{WIREs} Data Mining and Knowledge Discovery}, vol.~5, no.~5, pp.
  234--245, 2015.

\bibitem{2008Pareto}
Y.~Jin and B.~Sendhoff, ``Pareto-based multiobjective machine learning: An
  overview and case studies,'' \emph{IEEE Transactions on Systems Man \&
  Cybernetics Part C}, vol.~38, no.~3, pp. 397--415, 2008.

\bibitem{liu2021NC}
J.~Liu and Y.~Jin, ``Multi-objective search of robust neural architectures
  against multiple types of adversarial attacks,'' \emph{Neurocomputing}, vol.
  453, pp. 73--84, 2021.

\bibitem{liu2022TEVC}
J.~Liu, R.~Cheng, and Y.~Jin, ``Bi-fidelity evolutionary multiobjective search
  for adversarially robust deep neural architectures,''
  \emph{arXiv:2207.05321v1}, 2022.

\bibitem{wissam2015cyb}
W.~A. Albukhanajer, J.~A. Briffa, and Y.~Jin, ``Evolutionary multi-objective
  image feature extraction in the presence of noise,'' \emph{IEEE Transactions
  on Cybernetics}, vol.~45, no.~9, pp. 1757--1768, 2015.

\bibitem{EvoML-Survey}
\BIBentryALTinterwordspacing
A.~Telikani, A.~Tahmassebi, W.~Banzhaf, and A.~H. Gandomi, ``Evolutionary
  machine learning: A survey,'' \emph{ACM Comput. Surv.}, vol.~54, no.~8, oct
  2021. [Online]. Available: \url{https://doi.org/10.1145/3467477}
\BIBentrySTDinterwordspacing

\bibitem{ZHAN202242}
\BIBentryALTinterwordspacing
Z.-H. Zhan, J.-Y. Li, and J.~Zhang, ``Evolutionary deep learning: A survey,''
  \emph{Neurocomputing}, vol. 483, pp. 42--58, 2022. [Online]. Available:
  \url{https://www.sciencedirect.com/science/article/pii/S0925231222001345}
\BIBentrySTDinterwordspacing

\bibitem{ustun2019fairness}
B.~Ustun, Y.~Liu, and D.~Parkes, ``Fairness without harm: Decoupled classifiers
  with preference guarantees,'' in \emph{International Conference on Machine
  Learning}.\hskip 1em plus 0.5em minus 0.4em\relax PMLR, 2019, pp. 6373--6382.

\bibitem{8314455}
H.~Ishibuchi, R.~Imada, Y.~Setoguchi, and Y.~Nojima, ``Reference point
  specification in inverted generational distance for triangular linear pareto
  front,'' \emph{IEEE Transactions on Evolutionary Computation}, vol.~22,
  no.~6, pp. 961--975, 2018.

\bibitem{wang2021transfer}
X.~Wang, Y.~Jin, S.~Schmitt, and M.~Olhofer, ``Transfer learning based
  co-surrogate assisted evolutionary bi-objective optimization for objectives
  with non-uniform evaluation times,'' \emph{Evolutionary Computation}, pp.
  1--27, 2021.

\bibitem{wang2021transfer2}
X.~Wang, Y.~Jin, S.~Schmitt, M.~Olhofer, and R.~Allmendinger, ``Transfer
  learning based surrogate assisted evolutionary bi-objective optimization for
  objectives with different evaluation times,'' \emph{Knowledge-Based Systems},
  vol. 227, p. 107190, 2021.

\bibitem{7088618}
A.~Zhou and Q.~Zhang, ``Are all the subproblems equally important? resource
  allocation in decomposition-based multiobjective evolutionary algorithms,''
  \emph{IEEE Transactions on Evolutionary Computation}, vol.~20, no.~1, pp.
  52--64, 2016.

\bibitem{9541160}
A.~Chhabra, K.~Masalkovaitė, and P.~Mohapatra, ``An overview of fairness in
  clustering,'' \emph{IEEE Access}, vol.~9, pp. 130\,698--130\,720, 2021.

\bibitem{du2020fairness}
M.~Du, F.~Yang, N.~Zou, and X.~Hu, ``Fairness in deep learning: A computational
  perspective,'' \emph{IEEE Intelligent Systems}, vol.~36, no.~4, pp. 25--34,
  2020.

\bibitem{fitzsimons2019general}
\BIBentryALTinterwordspacing
J.~Fitzsimons, A.~Al~Ali, M.~Osborne, and S.~Roberts, ``A general framework for
  fair regression,'' \emph{Entropy}, vol.~21, no.~8, 2019. [Online]. Available:
  \url{https://www.mdpi.com/1099-4300/21/8/741}
\BIBentrySTDinterwordspacing

\bibitem{MatthewNips2016}
M.~Joseph, M.~Kearns, J.~H. Morgenstern, and A.~Roth, ``Fairness in learning:
  Classic and contextual bandits,'' in \emph{Proceedings of the 30th
  International Conference on Neural Information Processing Systems}, ser.
  NIPS'16.\hskip 1em plus 0.5em minus 0.4em\relax Red Hook, NY, USA: Curran
  Associates Inc., 2016, pp. 325--333.

\bibitem{singh2018fairness}
A.~Singh and T.~Joachims, ``Fairness of exposure in rankings,'' in
  \emph{Proceedings of the 24th ACM SIGKDD International Conference on
  Knowledge Discovery \& Data Mining}.\hskip 1em plus 0.5em minus 0.4em\relax
  New York, NY, USA: Association for Computing Machinery, 2018, pp. 2219--2228.

\bibitem{zehlike2017fa}
M.~Zehlike, F.~Bonchi, C.~Castillo, S.~Hajian, M.~Megahed, and R.~Baeza-Yates,
  ``Fa* ir: A fair top-k ranking algorithm,'' in \emph{Proceedings of the 2017
  ACM on Conference on Information and Knowledge Management}.\hskip 1em plus
  0.5em minus 0.4em\relax New York, NY, USA: Association for Computing
  Machinery, 2017, pp. 1569--1578.

\bibitem{komiyama2018nonconvex}
J.~Komiyama, A.~Takeda, J.~Honda, and H.~Shimao, ``Nonconvex optimization for
  regression with fairness constraints,'' in \emph{Proceedings of the 35th
  International Conference on Machine Learning}.\hskip 1em plus 0.5em minus
  0.4em\relax PMLR, 2018, pp. 2737--2746.

\bibitem{schmucker2020multi}
R.~Schmucker, M.~Donini, V.~Perrone, M.~B. Zafar, and C.~Archambeau,
  ``Multi-objective multi-fidelity hyperparameter optimization with application
  to fairness,'' in \emph{NeurIPS Workshop on Meta-Learning}, 2020.

\bibitem{padh2021addressing}
K.~Padh, D.~Antognini, E.~Lejal-Glaude, B.~Faltings, and C.~Musat, ``Addressing
  fairness in classification with a model-agnostic multi-objective algorithm,''
  in \emph{Uncertainty in Artificial Intelligence}.\hskip 1em plus 0.5em minus
  0.4em\relax PMLR, 2021, pp. 600--609.

\bibitem{martinez2020minimax}
N.~Martinez, M.~Bertran, and G.~Sapiro, ``Minimax {P}areto fairness: A multi
  objective perspective,'' in \emph{International Conference on Machine
  Learning}.\hskip 1em plus 0.5em minus 0.4em\relax PMLR, 2020, pp. 6755--6764.

\end{thebibliography}

%








\end{document}